\documentclass[onecolumn,10pt]{article}

\usepackage[utf8]{inputenc}
\usepackage[T1]{fontenc}
\usepackage{palatino}
\usepackage{textcomp}

\usepackage{amsmath}
\usepackage{amssymb}
\usepackage{mathtools}
\usepackage{gensymb}
\usepackage{siunitx}
\usepackage{cancel}
\usepackage{mhchem}

\usepackage{xcolor}
\usepackage{fullpage}
\usepackage{color,soul}
\usepackage{setspace}
\usepackage{lineno}
\usepackage{changepage}
\usepackage{pdflscape}
\usepackage{placeins}
\usepackage{fancyhdr}
\usepackage[bottom]{footmisc}
\usepackage[font=footnotesize]{caption}

\usepackage{graphicx}
\usepackage{subcaption}
\usepackage{booktabs}
\usepackage{array}
\usepackage{tabularx}
\usepackage{multirow}
\usepackage{adjustbox}
\usepackage{rotating}
\usepackage{wrapfig}

\usepackage{algorithm}
\usepackage{algpseudocode}

\usepackage{titling}
\usepackage{authblk}
\usepackage{subfiles}
\usepackage{chngcntr}
\usepackage{comment}

\usepackage{url}
\usepackage[square,sort,comma,numbers]{natbib}
\usepackage{hyperref}

\usepackage{lipsum}

\newcommand{\keywords}[1]{\par\vspace{0.5em}{\centering \textit{\textbf{Keywords: }}#1\par}}

\newcolumntype{L}[1]{>{\raggedright\let\newline\\\arraybackslash\hspace{0pt}}m{#1}}
\newcolumntype{C}[1]{>{\centering\let\newline\\\arraybackslash\hspace{0pt}}m{#1}}
\newcolumntype{R}[1]{>{\raggedleft\let\newline\\\arraybackslash\hspace{0pt}}m{#1}}
\title{\bfseries\fontsize{18}{20}\selectfont The \textit{Fuel Optimizer}: A Data-Driven Numerical Framework for Formulation of  Aviation Turbine Fuel}

\author[1,2]{Ana Larrañaga$^*$}
\author[1,2]{Steven L. Brunton}
\author[3]{Jacobo Porteiro}
\author[4]{Dario Lopez-Pintor}

\affil[1]{\small Department of Mechanical Engineering, University of Washington, Seattle, WA 98195, United States}
\affil[2]{\small NSF AI Institute in Dynamic Systems, University of Washington, Seattle, WA 98195, United States}
\affil[3]{\small CINTECX, Universidade de Vigo, Grupo de Tecnoloxía Enerxética (GTE), Vigo, 36310, Spain}
\affil[4]{\small Sandia National Laboratories, 7011 East Ave, Livermore, 94550, CA, United States}
\affil[*]{\small corresponding author: alarra@uw.edu}

\begin{document}

\date{}
\maketitle
\vspace{-.75in}
\begin{abstract}
\normalsize
The \textit{Fuel Optimizer} is an inverse design framework for sustainable aviation fuels that starts from a user-defined merit function and identifies the optimal combination of chemical species or hydrocarbon families that optimize a combination of targets. As a case study, a database of fuel blends meeting selected property standards was simulated in a reactor model to obtain pollutant emissions at cruise conditions. A surrogate model was developed to reduce the computational cost of evaluating candidate blends, taking fuel composition as input and predicting emissions as output. A genetic algorithm was used to optimize fuel formulation according to two merit functions designed to break the nitrogen oxides - CO trade-off, and minimize pollutant emissions over the landing-and-take-off cycle. Constraints included selected property standards (including seal swelling) and composition limits. Fuel candidates from the framework outperformed the training database across all merit functions, and the optimal candidates were validated through reactor simulations.
\end{abstract}
\vspace{0.2cm}

% ======================================================================
The aviation industry is actively working to diversify its fuel supply by increasing the use of synthetic aviation turbine fuels (SATFs), "drop-in" replacements for conventional jet fuel producible from renewable feedstocks~\cite{blakeyPROCI2011}. To meet the demands of commercial aviation, SATFs must comply with ASTM D1655~\cite{noauthor_standard_nodate} and ASTM D7566~\cite{noauthor_standard_d7566}, which define property requirements for turbine fuels in terms of safety, materials compatibility, and performance. SATFs must also contain specific hydrocarbon classes, including normal and iso-paraffins, aromatics, and cycloalkanes, to meet operational requirements, though their relative proportions directly affect engine emissions and performance~\cite{Maerkl2024ACP, Voigt2021ACP}.

To date, only eight SATF pathways have been approved under ASTM D7566. New pathways must follow the rigorous qualification process outlined in ASTM D4054~\cite{noauthor_standard_d4054}, which includes several rounds of property, material compatibility, and engine tests, estimated to cost around \$100M per candidate~\cite{calderonNREL2024}. Each new pathway is therefore a high-risk investment, and fuel designers must rely on predictive models to evaluate whether a candidate blend meets specifications before committing to testing-scale production. Existing models, however, lack the predictive capability to capture fuel effects on engine performance, partly owing to the complexity of fuel-engine interactions, for which machine learning (ML) approaches are particularly well-suited. 

ML techniques have emerged as powerful tools for modeling the complex, nonlinear interactions between fuel composition, combustion chemistry, and engine operating conditions, enabling more accurate prediction of NO$_{x}$, CO, unburned hydrocarbons, and soot than traditional semi-empirical approaches~\cite{zhouENERGY2022}. Larrañaga et al.~\cite{larranaga_data-driven_2026} developed a ML surrogate to predict the derived cetane number of complex fuel blends, capturing the interactions between physical properties and autoignition chemistry that govern ignition quality. Shao et al.~\cite{shaoCHEM2025} benchmarked multiple ML algorithms against density, net heat of combustion, viscosity, freeze point, and flash point, finding that the best-performing algorithm varied by property, reflecting differing degrees of non-linearity in the composition-property relationships. Aygun et al.~\cite{aygunENERGY2023} trained a convolutional neural network on operational and emissions data from a mixed-flow turbofan to predict fuel consumption, CO, and NO$_{x}$ as a function of engine pressure ratio and thrust, though fuel composition effects were not considered. Ma et al.~\cite{maTRANS2025} proposed a ML approach for real-time emission estimation across all flight stages using aircraft operational and atmospheric data as inputs, reporting emissions 22\% to 78\% higher than International Civil Aviation Organization (ICAO) reference values and identifying operational strategies with potential reductions of 5\% to 65\%; however, the study was limited to conventional kerosene at a single airport, leaving its generalizability unclear. To the best of the authors' knowledge, no numerical framework to date combines fuel design flexibility, sufficiently low computational cost, and acceptable accuracy in a single tool.

To address this gap, we developed the \textit{Fuel Optimizer}, a data-driven framework for the inverse design of SATF blends based on a user-defined merit function that optimizes pollutant emissions, fuel properties, or composition targets subject to user-defined constraints. The framework is based on a library of pure components including paraffins, iso-paraffins, cycloalkanes, aromatics, olefins, alcohols, ketones, furans, ethers, and esters. Engine-out emissions were generated via 1D reactor-network simulations calibrated against detailed combustor models and used to train a neural-network surrogate that evaluates candidate blends at a fraction of the computational cost. An optimization algorithm then explores the feasible composition space to identify optimum blends for a given merit function, with the goal of narrowing down promising candidates before committing to costly qualification testing.

%%%%%%%%%%%%%%%%%%%%%%%%%%%%%%%%%%%%%%%%%%%%%%%%%%%%%%%%%%%%
\vspace{-0.5cm}
\section{Background}
\label{sec:methods}
\vspace{-0.4cm}
\subsection{High-Fidelity Combustion and Emissions Modeling}
\label{subsec:hifi}
\vspace{-0.4cm}
A 1D reactor network model of a CFM56-7B27 turbofan engine was used to evaluate the performance and emissions of the multi-component fuel blends of the database~\cite{lopez-pintor_fuel_2026, mohammed_modeling_2026}. The reactor model was developed in ANSYS CHEMKIN and uses a series of perfectly stirred reactors and plug flow reactors interconnected to simulate the combustor, turbine and nozzle of the engine. The combustor is divided in six different zones (or reactors) that includes a rich-burn primary zone, a stoichiometric flame front, an intermediate zone and a dilution zone together with two wall zones to capture combustion inefficiencies and wall-quenching effects. The fuel and air flows into each zone were adjusted to match the fuel equivalence ratio in the primary zone and the mass-averaged temperature in the intermediate, dilution and wall zones from numerical simulations~\cite{zhang_predictions_2020}. The volume and residence time of each zone were estimated based on the adjusted flows and the geometric parameters of the engine~\cite{saboohi_development_2016}. The pressure ratios and geometric parameters of the turbine and nozzle were used to solve the temperature and pressure evolution within the turbomachinery, which was imposed in a series of plug flow reactors in CHEMKIN for prediction of engine nozzle-out emissions. 

Combustion chemistry was solved using a detailed chemical kinetic mechanism for aviation fuels. The gas-phase mechanism was developed by Lawrence Livermore National Laboratory (LLNL), consists of 8478 species and 33318 reactions, and is a combination of a detailed model for diesel- and aviation-like fuels from Wang et al.~\cite{wang_autoignition_2020} coupled with a polycyclic aromatic hydrocarbon (PAH) model from Kukkadapu et al.~\cite{kukkadapu_identification_2021} and a nitrogen oxides (NO$_x$) model from Glarborg et al.~\cite{glarborg_modeling_2018}. The mechanism has been validated against ignition delay time, laminar flame speed and reactor speciation measurements elsewhere~\cite{fang_fuel_2020, guzman_experimental_2019, richter_combined_2022}. Soot was modeled following a sectional method approach. This mechanism includes nucleation from gas-phase cyclo-penta pyrene (C18H10)~\cite{chung_computational_2011}, hydrogen abstraction carbon addition (HACA) growth from Kazakov et al.~\cite{kazakov_detailed_1995}, PAH condensation~\cite{veshkini_application_2016}, oxidation from Leistner et al.~\cite{leistner_kinetic_2012} and free-molecular coagulation. 

The model has been validated against engine emission measurements with Jet A from the International Civil Aviation Organization (ICAO) databank~\cite{noauthor_icao_nodate}, including NO$_x$, carbon monoxide (CO), particulate mass (PM) and particle number (PN), and a comparison between the experiments and simulations is included in the Supplementary Material. More details about this model and its performance can be found in~\cite{lopez-pintor_fuel_2026, garcia-oliver_development_2023, mohammed_modeling_2026}.

The 2508 multi-component blends of the fuel database and a surrogate fuel for Jet A proposed by Lopez-Pintor et al.~\cite{lopez-pintor_fuel_2026} were simulated in the reactor model at cruise-representative conditions (30$\%$ thrust, combustor inlet temperature of 613 K and combustor pressure of 11.3 bar). The exhaust temperature, the engine-out mole fractions of NO$_x$ (nitric oxide plus nitrogen dioxide), CO, carbon dioxide (CO$_2$) and water (H$_2$O), the soot mass flow rate, and the particle number density at the exhaust were incorporated into the database, which is included in the supplementary material of this paper. 

%%%%%%%%%%%%%%%%%%%%%%%%%%%%%%%%%%
\vspace{-0.4cm}
\subsection{Exploration vs Exploitation}
\label{subsec:optimization_background}
\vspace{-0.4cm}
Optimization tasks usually have to balance exploitation (i.e., the refinement of solutions within known high-performing regions) and exploration (i.e., broad sampling of the search space to avoid convergence to local optima)~\cite{goldberg_genetic_1989}. Gradient-based and deterministic methods are appropriate for exploitation but are ill-suited to high-dimensional, multimodal, or combinatorially complex scenarios, where the global optimum may be separated from any locally accessible solution. Genetic algorithms (GAs) address this limitation through a population-based search mechanism in which selection, crossover, and mutation operators act collectively to maintain diversity and cover a broad range of the solution space across successive generations~\cite{goldberg_genetic_1989}.

\begin{algorithm}
\caption{\small General genetic algorithm structure}
\label{alg:ga}
\begin{algorithmic}[1]
\small

\Require Population size $N$, generations $G$, fitness $f$, 
         feasible set $\mathcal{S}$, rates $p_c$, $p_m$
\Ensure  Best solution $x^*$

\State Initialize $P_0 = \{x_i\}_{i=1}^{N}$, where $x_i \in \mathcal{S}$
\For{$g = 1$ \textbf{to} $G$}
    \State $\varphi_i \leftarrow f(x_i)$ \quad $\forall\, x_i \in P_g$
    \State $P_g^{\text{sel}} \leftarrow \textsc{Select}(P_g, \{\varphi_i\})$
    \State $P_g^{\text{cross}} \leftarrow \textsc{Crossover}(P_g^{\text{sel}}, p_c)$
    \State $P_{g+1} \leftarrow \textsc{Mutate}(P_g^{\text{cross}}, p_m)$
    \State Discard $x_i \notin \mathcal{S}$ from $P_{g+1}$
\EndFor
\State $x^* \leftarrow \arg\min_{x_i \in \bigcup_g P_g} f(x_i)$
\State \Return $x^*$

\end{algorithmic}
\end{algorithm}

The population-based nature of GAs (see Algorithm~\ref{alg:ga}) retains information not only about the best solutions, but about the full distribution of evaluated candidates. In the present work, the chemical space of possible blends is high-dimensional and combinatorially large, making a complete evaluation intractable. The main objective is not to converge to a single optimal blend, but to identify which chemical families and pure fuel components persistently appear in high-performing formulations. By analyzing candidate composition across generations, it is possible to infer which regions of chemical space are recurrently selected by the fitness criterion, orienting the search toward the most promising blend components~\cite{mueller_optimization_2022}.

%%%%%%%%%%%%%%%%%%%%%%%%%%%%%%%%%%%%%%%%%%%%%%%%%%%%%%%%%%%%
\vspace{-0.4cm}
\section{Methodology}
\label{sec:methodology}
\vspace{-0.4cm}
\subsection{Fuel Component Library and Blend Generation}\label{library}
\vspace{-0.4cm}
SATF candidates were formulated by combining 51 single component species that included n-alkanes (carbon chain length from 9 to 20), iso-alkanes (single and multi-branched), cyclo-alkanes (mono and di-cyclic), aromatics (mono and di-cyclic), olefins (mono and dienes), alcohols (normal and branched), ethers, esters, ketones and furans. These 51 species were selected to cover a wide range of molecular structures while considering the availability of chemical kinetic models and their ability to meet ASTM D1655 constrains~\cite{noauthor_standard_nodate}. A brute-force Monte Carlo approach was implemented to formulate 2508 multi-component blends with 3 to 5 components and with a minimum species content of 2$\%_{mole}$ that met the ASTM D1655 requirements for density, viscosity, distillation, freeze point, flash point and neat heat of combustion. To minimize database bias, all the species of the palette were forced to be present in at least 50 of the multi-component blends. More details of the fuel database and fuel property estimation methods are included in Supplementary Material.   

%%%%%%%%%%%%%%%%%%%%%%%%%%%%%%%%%%%%%%%%%%%%%%%%%%%%%%%%%%%%
\vspace{-0.4cm}
\subsection{Surrogate Modeling of 1D Reactor Plume Model}
\label{subsec:ml}
\begin{figure}[t]
    \centering
    \vspace{-4cm}
    \includegraphics[trim={1cm 0 1cm 0},clip,width=1\linewidth]{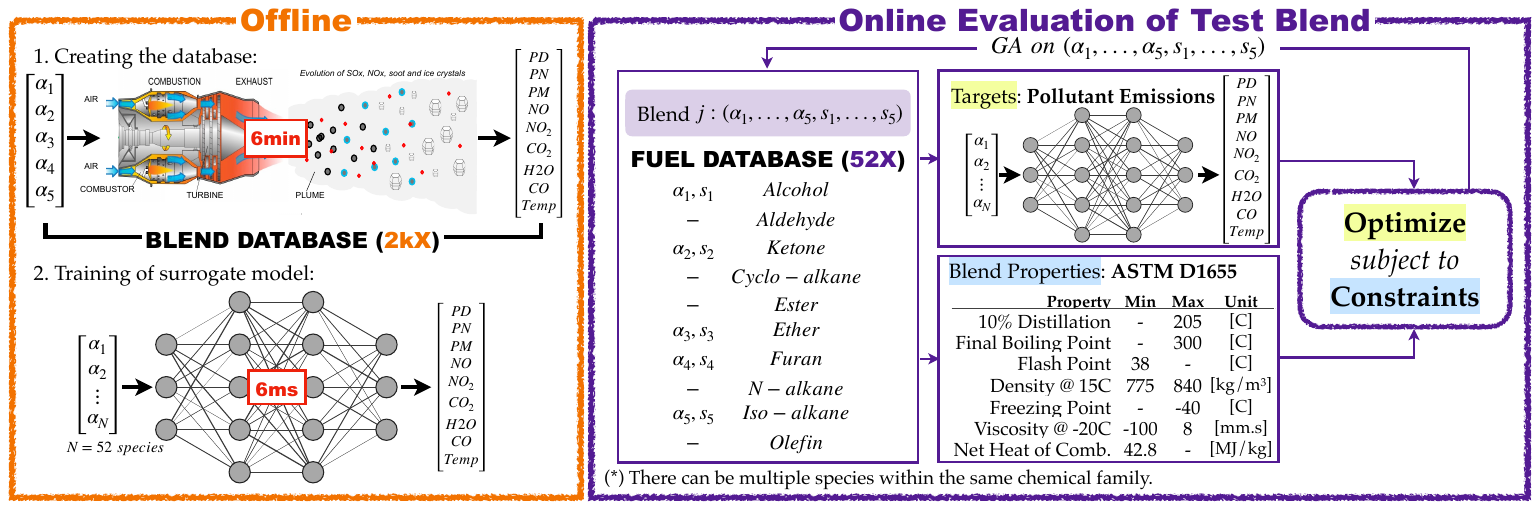}
    \scriptsize
    \vspace{-4.5cm}
    \caption{This work consists of two separate tasks. The first one corresponds to the offline training of the surrogate model (orange), which includes the generation of a fuel blend database containing all chemical species in the pure components database and the subsequent training of the neural network. The second stage corresponds to the online optimization procedure, in which the trained surrogate model is embedded within a genetic algorithm loop. Each candidate blend is evaluated using the surrogate model to compute the corresponding merit functions, MF1 and MF2, while the blends are constrained according to the Standard Specification for Aviation Turbine Fuels (ASTM D1655). All candidate blends evaluated by the algorithm, together with their corresponding merit function values, are stored and tracked to analyze which chemical families are preferentially selected for each optimization objective.}
    \vspace{-0.5cm}
    \label{workflow}
\end{figure}

\vspace{-0.4cm}
One of the main limitations of the original 1D reaction plume model for optimization tasks is its computational cost, as a single evaluation requires several minutes on a standard laptop. Since genetic algorithms typically evaluate thousands of candidate blends per iteration, direct integration of the model into the optimization framework becomes computationally impractical. To overcome this limitation, a surrogate model was developed to predict the same outputs as the original model from the proportions of the chemical species in each blend (see Fig.~\ref{workflow}).

Given to the nonlinear nature of the underlying physics and the available dataset, a neural network architecture was adopted. Because the optimization objective is primarily based on emission-related metrics used in the merit function, special attention was given to the prediciton accuracy for PN, PM, CO and NOx, rather than for the complete set of outputs. Relative to the original 1D model, the surrogate reduces evaluation time by approximately three orders of magnitude. Additional details on the model architecture and hyperparameter tuning are provided in the Supplementary Material.

%%%%%%%%%%%%%%%%%%%%%%%%%%%%%%%%%%%%%%%%%%%%%%%%%%%%%%%%%%%%
\vspace{-0.4cm}
\subsection{Optimization}
\label{subsec:opt}
\vspace{-0.4cm}
As previously mentioned, a genetic algorithm (GA) was selected to encourage exploration of the feasible fuel-blend space. The optimization vector consisted of 10 elements: five corresponding to the species identifiers in the database and five corresponding to the mass fractions of each species in the blend. A minimum mass fraction of 2\% was imposed to ensure that each component had a measurable effect on the objective function, while the sum of all mass fractions was constrained to unity.

The GA hyperparameters are reported in the Supplementary Material. In practice, the algorithm was allowed to explore the widest possible range of blends for each species in the database. For this purpose, one species identifier was fixed to guarantee its presence in the blend, whereas the remaining nine parameters were optimized. As discussed later in this section, additional constraints were introduced to ensure realistic fuel properties and practically feasible blends.

It is important to note that all candidate blends evaluated by the GA are stored and tracked to enable an in-depth analysis of the chemical families most frequently selected for the optimization of the different merit functions. In a future research line, this information could be used to introduce an additional optimization stage focused on the best-performing species and associated chemical families, thereby reducing the search space. 

%%%%%%%%%%%%%%%%%%%%%%%%%%%%%%%%%%%%%%%%%%%%%%%%%%%%%%%%%%%
\vspace{-0.4cm}
\subsubsection{Merit Functions}
\label{sec:mf}
\vspace{-0.4cm}
The \textit{Fuel Optimizer} has been tested using a single objective that combines relevant targets predicted by the surrogate model. The analysis focuses on two different merit functions (MFs):

\begin{equation}
\label{eq:MF_general}
MF_k = \sum_{i} w_i \frac{i - i_{\text{JetA}}}{\sigma_{i,\text{DB}}}, \qquad k \in \{1,2\}
\end{equation}

where $\sigma_{i,\text{DB}}$ is the corresponding standard deviation derived from the fuel blend database, and $i_{\text{JetA}}$ denotes the baseline for target $i$ corresponding to Jet A. For both merit functions, negative values, which represent an improvement with respect to Jet A, are desired.  

\vspace{0.5cm}
\begin{itemize}
    \item \textbf{CO-NO$_x$ trade-off}: $MF_1$ is intended to break the CO-NO$_x$ emissions trade-off typically seen in jet engines, using equal weighting $w_{CO} = w_{NO_x} = 0.5$, such that
    \[
    MF_1 = w_{CO} \frac{CO - CO_{\text{JetA}}}{\sigma_{CO,\text{DB}}}
         + w_{NO_x} \frac{NO_x - NO_{x,\text{JetA}}}{\sigma_{NO_x,\text{DB}}},
    \]
    Variations in operating parameters that tend to reduce NO$_x$, such as decreasing the residence time or the combustor temperature, also lead to incomplete combustion raising CO emissions and vice versa. Even though this trade-off is not linear, for simplicity, the same weight was imposed to both CO and NO$_x$ in this case study.
    
    \item \textbf{LTO regulated emissions}: $MF_2$ is intended to minimize regulated emissions. Weighting factors, $w_{PN} = 0.12$, $w_{PM} = 0.16$, and $w_{NO_x} = 0.72$, have been defined based on how the emission values over the LTO cycle of the CFM56-7B27 turbofan engine reported in the ICAO databank compare against the legal limits. The targets are PN, PM, and NO$_x$, such that
    \[
    MF_2 = w_{PN} \frac{PN - PN_{\text{JetA}}}{\sigma_{PN,\text{DB}}}
         + w_{PM} \frac{PM - PM_{\text{JetA}}}{\sigma_{PM,\text{DB}}}
         + w_{NO_x} \frac{NO_x - NO_{x,\text{JetA}}}{\sigma_{NO_x,\text{DB}}},
    \]
    For this engine, NO$_x$ emission limits are significantly more difficult to meet than PN and PM emission limits, and this has been reflected in the weighting factors of $MF_2$.  
\end{itemize}

\begin{figure}[t]
    \centering
    \includegraphics[width=0.9\linewidth]{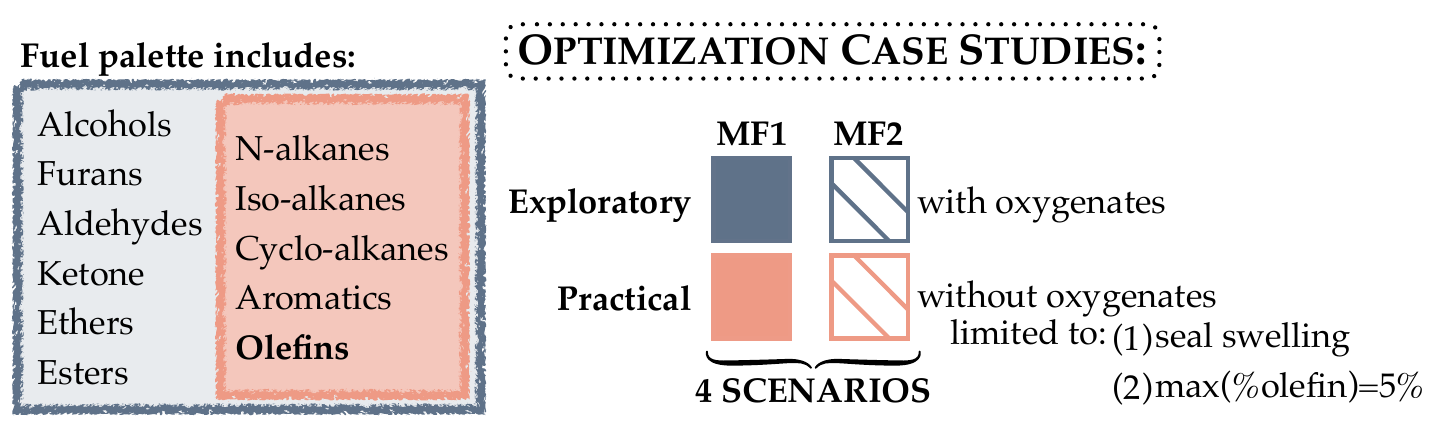}
    \scriptsize
    \caption{This study explores two fuel blend optimization scenarios: (1) fuel blends that contain oxygenates (purple), and (2) fuel blends without oxygenates but subject to additional constraints (green). The second case is a more practical optimization scenario while the first is more exploratory. Two merit functions, MF1 and MF2, are optimized. MF1 focuses on CO-NOx trade-off and MF2 focused on LTO regulated emissions.}
    \label{fig:scenarios}
    \vspace{-.4cm}
\end{figure}
%%%%%%%%%%%%%%%%%%%%%%%%%%%%%%%%%%%%%%%%%%%%%%%%%%%%%%%%%%%
\vspace{-.6cm}
\subsubsection{Constraints}
\label{sec:const}
\vspace{-.4cm}
The \textit{Fuel Optimizer} allows the user to define a set of constraints for the optimization algorithm. The number of components in each fuel candidate was set to 5, with a minimum mole content of 2$\%$ for each individual component. It is important to note that the palette of species available for the generation of blend candidates, and the number of components in each fuel candidate, are variable and can be modified depending on the interest on specific chemical species or families.

The GA input vector is structured such that the first five elements represent the molar proportions of the fuel components, while the remaining five elements specify the corresponding component indices in the database. The molar proportions are constrained to sum to 1.

Fuel candidates must also meet the ASTM D1655 requirements for density, viscosity, distillation, freeze point, flash point and neat heat of combustion (same physical property constraints used in the generation of the database, see section \ref{library} for details). In a more \textit{exploratory} approach, no other constraints were imposed and the \textit{Fuel Optimizer} was allowed to include species from all chemical families of the palette. Then, a more \textit{practical} approach was considered where oxygenated components were eliminated from the palette and the maximum olefinic content of fuel candidates was restricted to 5$\%$ as per ASTM D1655 limits. Additionally, nitrile seal swelling properties of the fuel components of the palette have been obtained from the literature~\cite{Faulhaber2023}\cite{Hamilton2023} and used to estimate the seal swell of fuel candidates using a linear blending rule based on liquid volume fractions~\cite{Faulhaber2023}. In this second approach, fuel candidates must show a minimum swell rate equal to that of conventional kerosene (defined as 12.4\% $\pm$ 2.3\% from 16 different samples from the literature~\cite{Faulhaber2023}\cite{Hamilton2023}), ensuring fuel system functionality. These two sets of constraints will be referred to as the \textit{exploratory} approach and the \textit{practical} approach hereafter, respectively (see Fig.~\ref{fig:scenarios}).
%%%%%%%%%%%%%%%%%%%%%%%%%%%%%%%%%%%%%%%%%%%%%%%%%%%%%%%%%%%%
\vspace{-.4cm}
\section{Results \& Discussion}
\label{sec:results}
\vspace{-.4cm}
This section presents the results obtained with the surrogate models developed for the different versions of the database, together with the optimization results. The GA results are analyzed to examine the diversity of the database and to validate the optimal blends, ensuring that their performance is consistent with the predicted behavior by the surrogate model.

%%%%%%%%%%%%%%%%%%%%%%%%%%%%%%%%%%%%%%%%%%%%%%%%%%%%%%%%%%%%
\vspace{-.4cm}
\subsection{Performance Evaluation}
\label{subsec:performance_evaluation}
\vspace{-.4cm}
Two surrogate models were trained independently for the two \textit{Fuel Optimizer} configurations. The models achieved prediction errors below 10\% for the targets relevant to MF1 and MF2, which include particulate mass (PM), particle number (PN), NO$_x$ and CO (see the Supplementary Material for details on model accuracy). Once trained, the models were included in the optimization loop, as  shown in Fig~\ref{workflow}. 

\begin{figure}[t]
    \centering
    \includegraphics[width=0.8\linewidth]{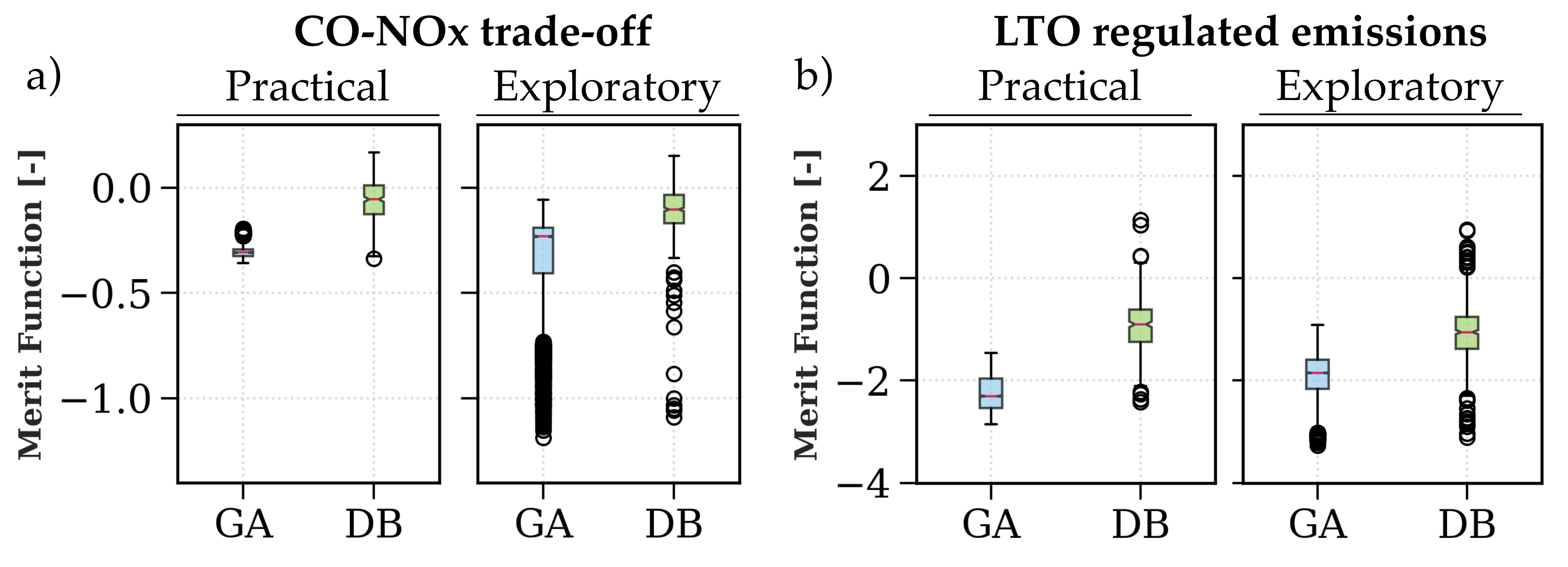}
    \scriptsize
    \vspace{-0.2cm}
    \caption{Comparison between the brute-force database (green, labeled DB) and the feasible candidate blends evaluated by the GA (blue, labeled GA) for all optimization scenarios. The left panels correspond to the CO--NO\textsubscript{x} trade-off (MF1), whereas the right panels correspond to the LTO regulated-emissions objective function (MF2). The bottom panels represent the practical approach, in which oxygenates are excluded and constraints on olefin content and seal-swelling performance are imposed, while the top panels correspond to the exploratory approach with fewer constraints. The mean merit-function value and standard deviation for each case are indicated in the corresponding panels.}
    \label{results-1}
    \vspace{-0.5cm}
\end{figure}

Fig.~\ref{results-1} shows a comparison of the performance of the fuel candidates obtained with the GA against those of the database obtained by brute force for the four cases described in Fig.~\ref{fig:scenarios}. Results highlight the ability of the GA to generate thousands of fuel candidates that are better than those obtained from a brute force approach, leading to merit function mean values that are systematically more negative (i.e., better) than those of the database (for the reader's reference, conventional kerosene shows a merit function value equal to 0). The GA improved the results of the database in all cases, with the optimum fuel from the GA leading to better merit function values than the optimum fuel of the database. 

%%%%%%%%%%%%%%%%%%%%%%%%%%%%%%%%%%%%%%%%%%%%%%%%%%%%%%%%%%%%
\vspace{-.4cm}
\subsection{Chemical Family Analysis for GA Results}
\begin{figure}[t]
    \centering
    \includegraphics[width=1\linewidth]{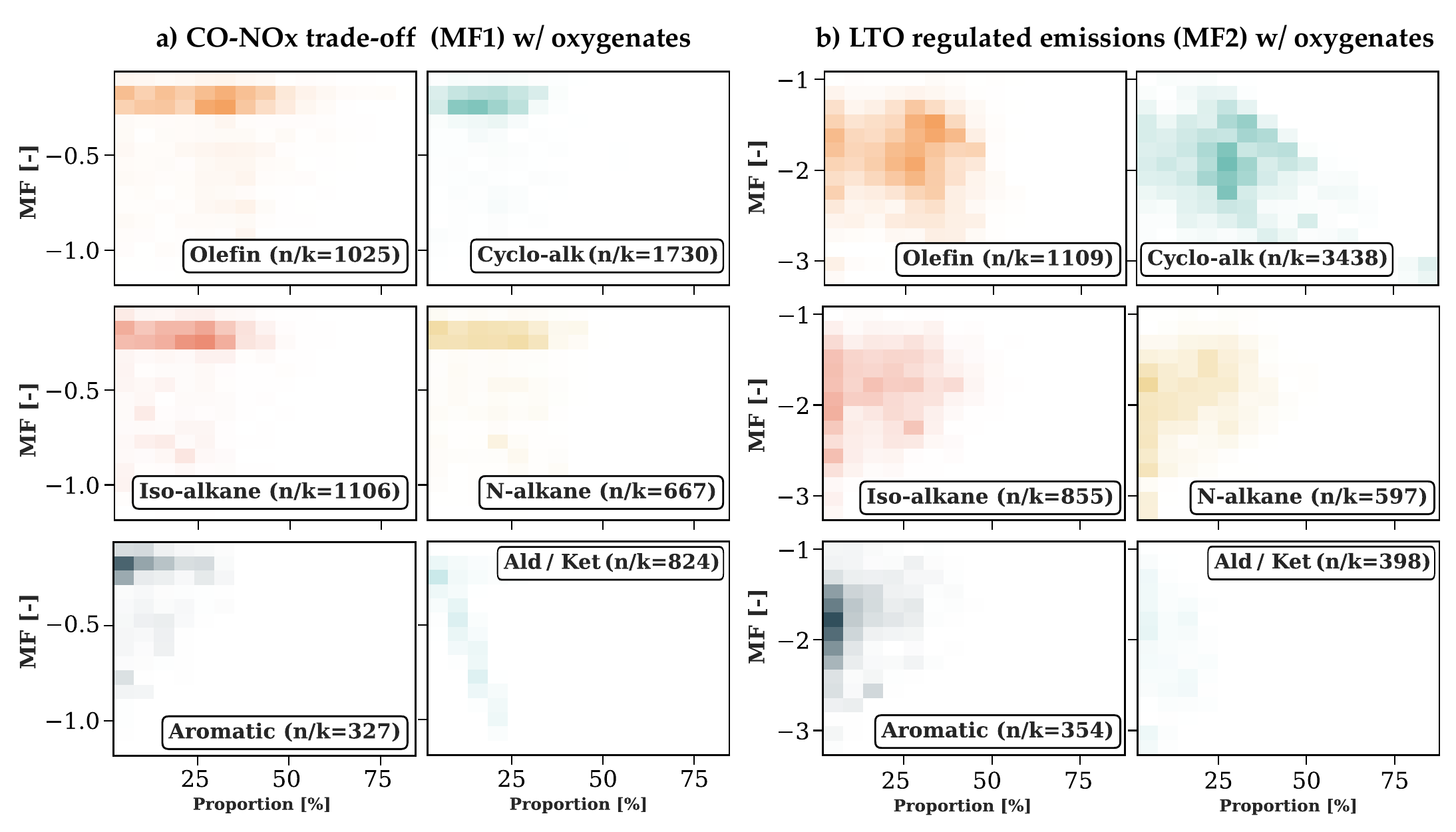}
    \scriptsize
    \vspace{-0.5cm}
    \caption{Results for the \textit{exploratory approach} analyzed by chemical family, focusing on the most relevant families included in the database (additional results are provided in the Supplementary Material). Figure (a) on the left corresponds to the CO--NO\textsubscript{x} trade-off optimization, whereas figure (b) on the right corresponds to the LTO regulated-emissions optimization. The merit-function value is plotted on the y-axis, while the proportion of each chemical family in the blend is plotted on the x-axis. The color intensity represents the frequency of occurrence, with darker regions indicating a higher density of blends within a given composition range. The annotations in each subplot indicate the ratio between the number of occurrences of a given family among all candidate blends ($n$) and the number of species of that family available in the database ($k$), which serves as an indicator of the preference of the GA for each family.}
    \vspace{-0.5cm}
    \label{results-3}
\end{figure}

\vspace{-.4cm}
\textbf{Exploratory Fuel Optimizer.} Fig.~\ref{results-3} shows the occurrence frequency of each family among all fuel candidates generated by the GA. Frequency is represented by color intensity in a merit function vs. proportion map for MF1 (left) and MF2 (right) for the exploratory approach. The number of times each family is present among all fuel candidates (n) divided by the number of species of each family available in the palette (k), which is a metric of the preference of the GA for each family, is included in the figure. Alcohols, furans, esters and ethers are significantly less selected than the other families, so they are not included in the figure (but results for these families are included in the supplementary material). Frequency results indicate that cyclo-alkanes are the most selected family for both merit functions (especially for MF2). This is because n-butyl cyclohexane acts as a super-component that meets all ASTM property constraints while having combustion properties that align well with the MF1 and MF2 targets (as it will be discussed in detail below). 
In general, all families were able to achieve similar minimum merit function values, indicating that \textit{good} fuel candidates can be formulated from any chemical family with the right combination of components (exceptions include furans, esters and ethers for MF1, and furans for MF2, see supplementary materials for details). However, Fig.~\ref{results-3} shows clear preferences in the proportion in which each chemical family is present among the best fuel candidates. 

For MF1 (NO$_x$-CO tradeoff), the best candidates show a clear preference for increasing the ketone content and for a low aromatic content combined with moderate amounts of iso-alkanes and olefins. In general, higher air-fuel stoichiometric ratios lead to more diluent per unit of chemical energy released, lower adiabatic flame temperatures and lower thermal NO$_x$ formation, which is associated with better MF1 values. Similarly, aviation fuels with lower neat heat of combustion also tend to show lower adiabatic flame temperatures and lower NO$_x$, with oxygen content being a dominant factor affecting this property. On the other hand, CO emissions are controlled by the ability of the fuel to burn in the wall regions and to complete the CO-to-CO$_2$ oxidation reactions within the combustor residence time. Fuels with higher flame speed generally show higher resistance to quenching and require less residence time to complete combustion. Therefore, the genetic algorithm was expected to maximize the air-fuel stoichiometric ratio and oxygen content of the fuel to minimize NO$_x$, and maximize flame speed to minimize CO emissions. 

Among all the species of the palette, cyclopentanone (ketone) stands out for its high oxygen content and flame speed, and has physical properties that allow blending at relatively high fractions ($\sim$20$\%$). Thus, for the exploratory approach where oxygenated species were considered, cyclopentanone was identified as the key species to break the NO$_x$-CO tradeoff (note how the MF1 value systematically improves as the ketone content increases in Fig.~\ref{results-3}). 

For MF2 (LTO emissions), the best candidates are achieved with very large concentrations of cyclo-alkanes and low aromatic and iso-alkane content. As expected, MF2 improved as the aromatic content and, to a less extent, the iso-alkane content of the fuel candidates decreased because of the negative impact of these chemical families on soot emissions. Cyclo-alkanes, and more specifically n-butyl cyclohexane, have been identified as the key chemical family to decrease regulated emissions, with MF2 decreasing as the n-butyl cyclohexane content increased. This is because the slower burn rate of cyclo-alkanes led to lower combustor temperatures and, therefore, lower NO$_x$. As it will be shown later, this slower burn rate of cyclo-alkanes led to higher CO emissions, but CO was not included in MF2 and, therefore, it was ignored by the optimization algorithm.

\begin{figure}[t]
    \centering
    \includegraphics[width=1\linewidth]{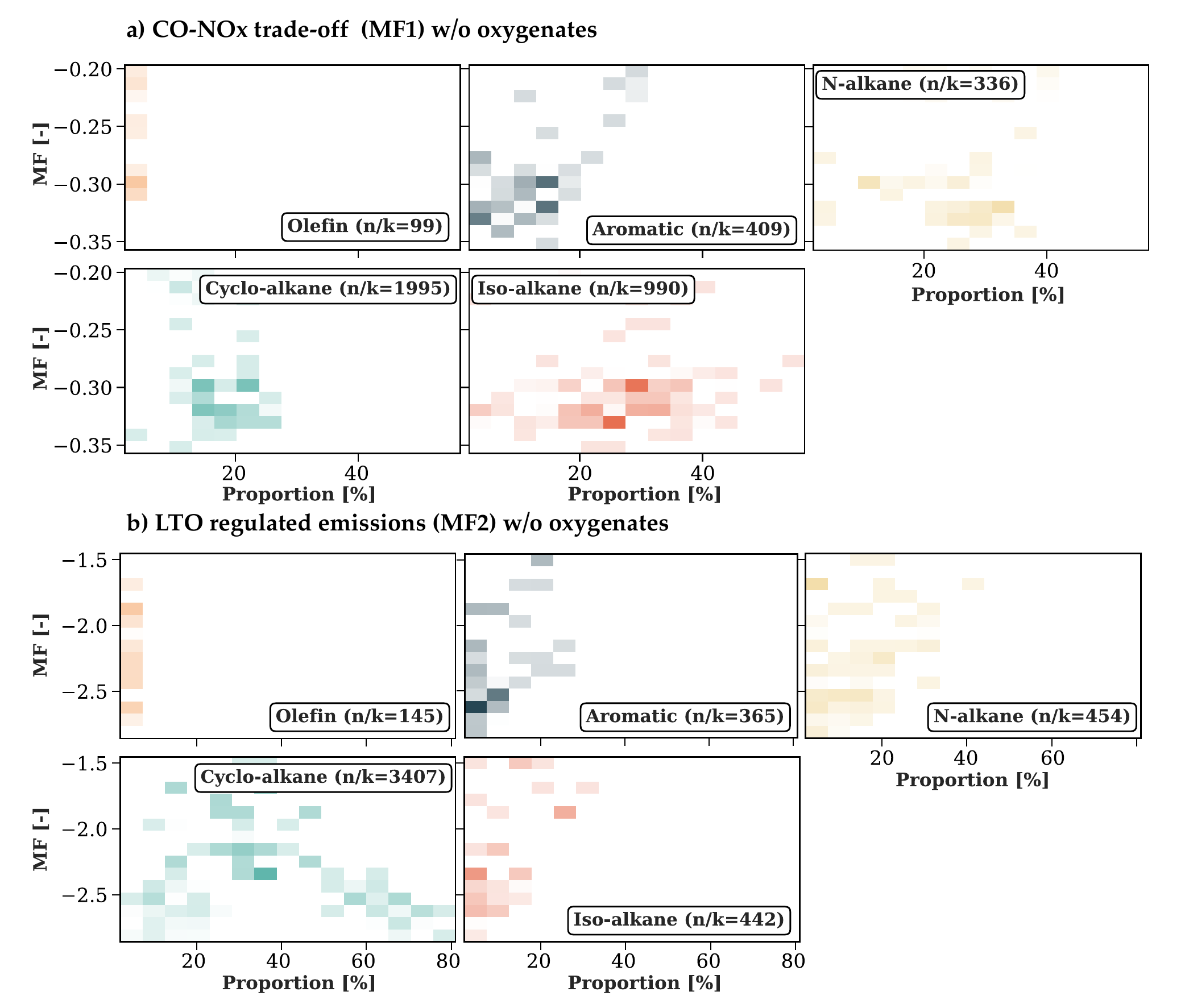}
    \scriptsize
    \caption{Results for the \textit{practical approach} analyzed by chemical family, considering that in this case the database is limited to not include oxygenates. Figure (a) on the top corresponds to the CO--NO\textsubscript{x} trade-off optimization, whereas figure (b) on the bottom corresponds to the LTO regulated-emissions optimization. The merit-function value is plotted on the y-axis, while the proportion of each chemical family in the blend is plotted on the x-axis. The color intensity represents the frequency of occurrence, with darker regions indicating a higher density of blends within a given composition range. The annotations in each subplot indicate the ratio between the number of occurrences of a given family among all candidate blends ($n$) and the number of species of that family available in the database ($k$), which serves as an indicator of the preference of the GA for each family.}
    \label{results-2}
    \vspace{-0.4cm}
\end{figure}

\textbf{Practical Fuel Optimizer.} Fig.~\ref{results-2} shows the occurrence frequency of each family among all fuel candidates generated by the GA with frequency represented by color intensity in a merit function vs. proportion map. Results for MF1 (top) and MF2 (bottom) are shown for the practical approach. Similarly to figure \ref{results-3}, the ratio between the number of times each family is present among all fuel candidates and the number of species of each family available in the palette (n/k) is included in the figure. Cyclo-alkanes are again the most selected family for both merit functions, with the GA showing a very strong preference for this family for the MF2 case. On the other hand, olefins are significantly less used than other chemical families, consistent with the limitation of maximum 5$\%$ olefinic content imposed in the practical approach. 
Except for the olefins, all the chemical families show similar mean and minimum merit function values, indicating that near-optimum fuel candidates are achieved by combining components from all chemical families. 

An analysis of the proportion of each chemical family against its merit function value indicate that, for MF1 (NO$_x$-CO tradeoff), there is not single chemical family that is strongly preferred in high concentrations, with the optimum blend resulting from the combination of aromatics, cyclo-alkanes, iso-alkanes and n-alkanes. (Note that ketones are not a blending option in the practical approach). MF1 tends to improve as the aromatic content of the candidate decreases because aromatics have low burning rates and tend to increase CO emissions.  

For MF2 (LTO emissions), the merit function value improves as the cyclo-alkane content increases up to concentrations around 80$\%$ and the content of the other chemical families decreases. As discussed before, n-butyl cyclohexane acts as a \textit{super component} that can be blended in very large proportions while meeting ASTM constraints and leading to lower NO$_x$ thanks to lower combustion temperatures due to its slower burn rate. MF2 also benefits from low aromatic and iso-alkane proportions due to the high soot propensity of these components (especially aromatics). Nevertheless, a minimum amount of aromatics is required to meet the seal swelling constraint with the best blends showing 4 - 8$\%$ aromatics combined with larger amounts of cyclo-alkanes. 

%%%%%%%%%%%%%%%%%%%%%%%%%%%%%%%%%%%%%%%%%%%%%%%%%%%%%%%%%%%%
\vspace{-0.2cm}
\subsection{Top Best Chemical Families and Species for GA Results}
\vspace{-0.4cm}
\begin{figure}[t]
    \centering
    \includegraphics[width=1\linewidth]{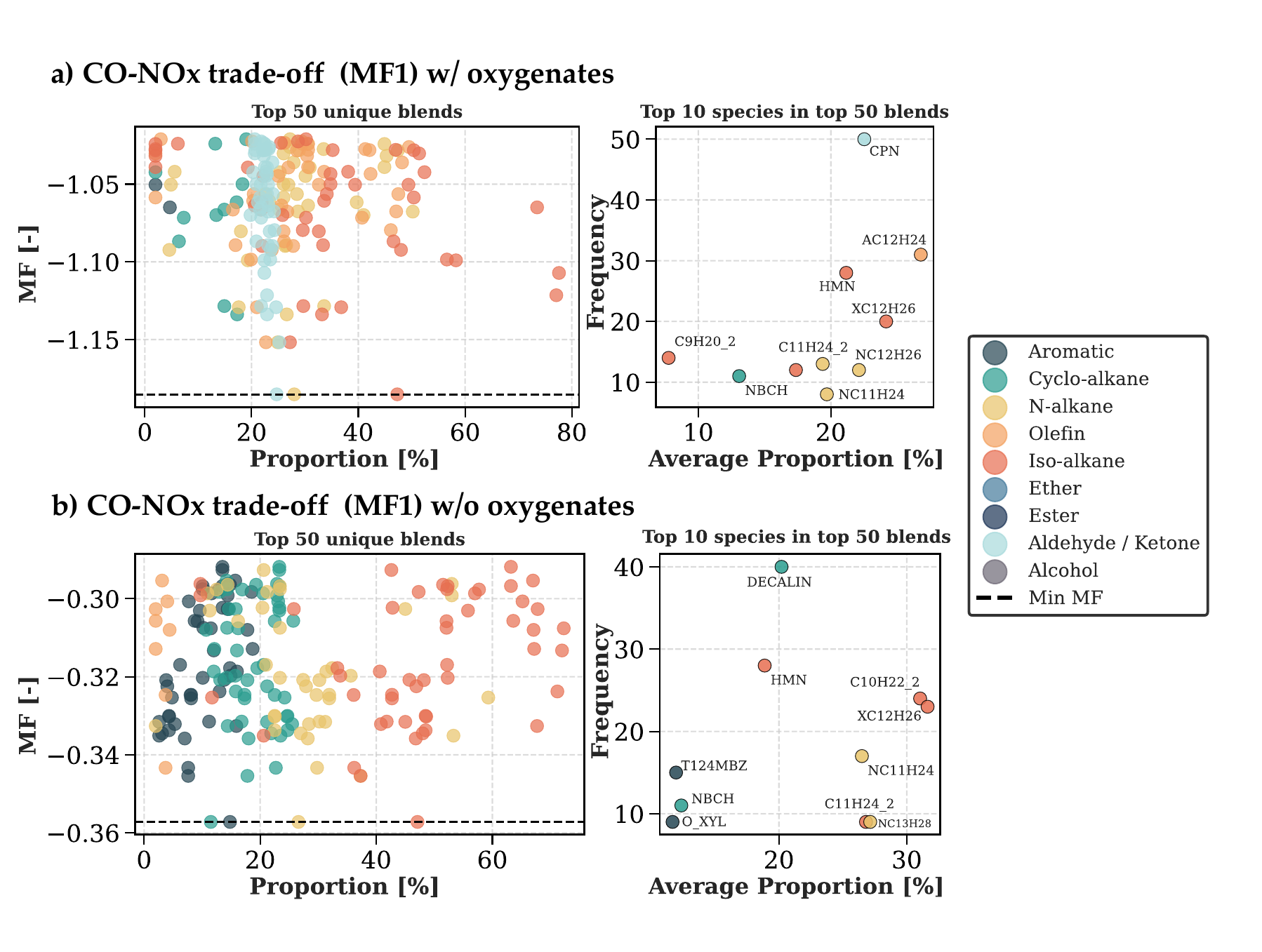}
    \scriptsize
    \vspace{-0.7cm}
    \caption{Detailed analysis of the 50 best unique fuel blends by chemical family (left panels) and by individual chemical species (right panels) for the CO--NO\textsubscript{x} trade-off objective function. The top panels correspond to the exploratory approach, which considers the complete database, whereas the bottom panels correspond to the practical approach, in which oxygenates are excluded and additional constraints on olefin content and seal-swelling performance are imposed. The minimum merit-function value is indicated in the left panels by a dashed black line, following the same analysis presented in Figs.~\ref{results-3} and ~\ref{results-2}  for the relationship between merit-function value and chemical-family composition. The right panels show the average contribution of each species to the top 10 blends together with the frequency with which each species appears among the optimal blends.}
    \label{results-4}
    \vspace{-0.4cm}
\end{figure}

The hydrocarbon family distribution of the top 50 best performing fuel candidates is shown in Fig.~\ref{results-4} against their merit function values for MF1 and for the exploratory (top) and practical (bottom) approaches. The figure shows also the occurrence frequency and proportion of the 10 most common species within the top 50 candidates. For the exploratory approach, cyclopentanone (CPN, a ketone) is present in all 50 best blends with a proportion close to 25$\%$. Larger CPN fractions were not viable due to its low energy density (32.6 MJ/kg), which would lead to fuel candidates that would not meet the ASTM limits for lower heating value (min. 42.8 MJ/kg). In these top performing blends, CPN is mixed with moderate amounts of large olefins and iso-alkanes (C12 - C16) that provide high air/fuel stoichiometric ratios (for low NO$_x$ emissions), while maintaining high burning rates (for low CO emissions) and the physical properties required to meet constraints. 

For the practical approach, decahydronaphthalene (decalin, a cyclo-alkane) is the is the most common species among the top performing candidates. Decalin is selected by the GA as a compromise alternative to aromatics to meet the seal swelling requirement while maintaining low CO emissions (note that aromatics tend to increase CO due to a very low reactivity and slow burning rate). Similarly to the exploratory approach, decalin is blended with moderate amounts of iso-alkanes or, less frequently, n-alkanes for high air/fuel stoichiometric ratios (for low NO$_x$ emissions) and to meet physical property constraints. Interestingly, low amounts of aromatics are sometimes selected, likely to meet seal swelling constraints. When this occurs, the GA selects highly branched mono-aromatics (trimethyl benzene or xylene) that show higher laminar flame speeds than non-branched mono- and di-aromatics.

\begin{figure}[t]
    \centering
    \includegraphics[width=1\linewidth]{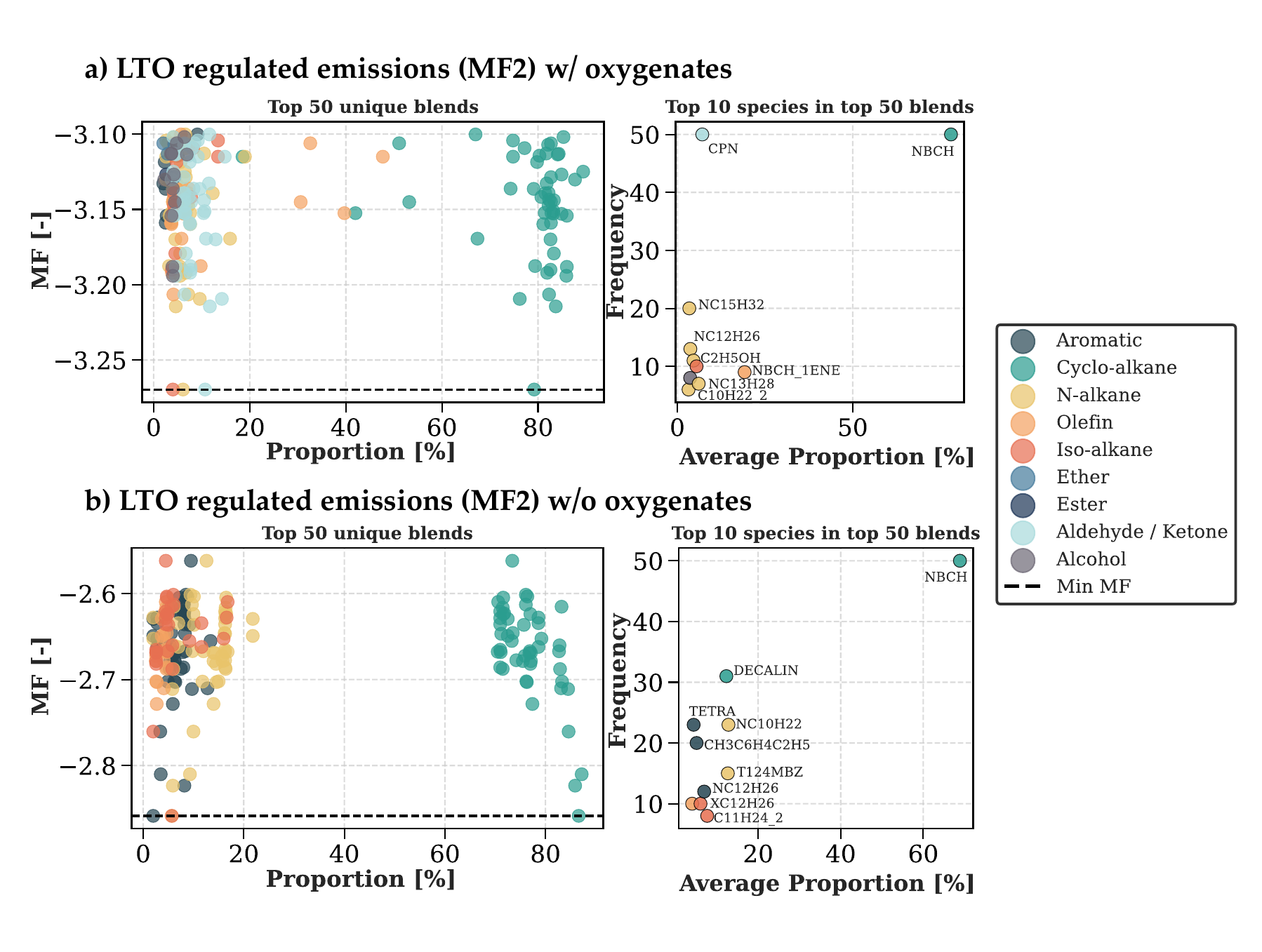}
    \scriptsize
    \vspace{-0.5cm}
    \caption{Detailed analysis of the 50 best unique fuel blends by chemical family (left panels) and by individual chemical species (right panels) for the LTO regulated emissions objective function. The top panels correspond to the exploratory approach, which considers the complete database, whereas the bottom panels correspond to the practical approach, in which oxygenates are excluded and additional constraints on olefin content and seal-swelling performance are imposed. The minimum merit-function value is indicated in the left panels by a dashed black line, following the same analysis presented in Figs.~\ref{results-3} and ~\ref{results-2}  for the relationship between merit-function value and chemical-family composition. The right panels show the average contribution of each species to the top 10 blends together with the frequency with which each species appears among the optimal blends.}
    \label{results-5}
    \vspace{-0.4cm}
\end{figure}

The hydrocarbon family distribution of the top 50 best performing fuel candidates for MF2 is shown in Fig.~\ref{results-5} against their merit function values for the exploratory (top) and practical (bottom) approaches. Similarly to Fig.~\ref{results-4}, the Fig.~includes the occurrence frequency and proportion of the 10 most common species within the top 50 candidates.
For MF2 and the exploratory approach, the top 50 best performing blends are all composed of very large fractions of n-butyl cyclohexane (NBCH, average proportion of 78$\%$) with mild amounts of CPN (average proportion of 7.1$\%$). NBCH and CPN both have low soot propensity and tend to decrease NO$_x$ due to either a slow burning rate (NBCH) or a high oxygen content (CPN). The GA tends to blend these components with multiple n-alkanes to \textit{fill the gaps} while keeping soot propensity low.  

For MF2 and the practical approach, NBCH is again a \textit{super component} present in very large fractions (approx. 70$\%$) that provides low soot propensity, lower NO$_x$, and physical properties within ASTM limits. In this case, large amounts of cyclo-alkanes are combined with a minimum amount of aromatics to meet the seal swelling constraint, with mono-aromatic species (methylethyl benzene and tetrahydronapththalene) being preferred over di-aromatics because they represent a better compromise between seal swelling and soot propensity.

%%%%%%%%%%%%%%%%%%%%%%%%%%%%%%%%%%%%%%%%%%%%%%%%%%%%%%%%%%%%
\vspace{-0.4cm}
\subsection{Validation of Optimum Blends}
\label{subsec:validation}

\begin{table}[t]
\centering
\caption{Detailed composition of optimum fuel blends for the CO--NO\textsubscript{x} trade-off (MF1) and LTO regulated emissions (MF2) in the different optimization case scenarios: (1) top section corresponds to the exploratory approach, and (2) bottom section corresponds to the practical approach. See the Supplementary Material for a glossary of species identifiers.}
\label{tab:fuel_blends}
\renewcommand{\arraystretch}{1}
\begin{adjustbox}{max width=0.9\textwidth}
\begin{tabular}{c l l c | l l c}
\toprule
\multicolumn{4}{c|}{\textbf{CO-NOx trade-off (MF1)}} &
\multicolumn{3}{c}{\textbf{LTO regulated emissions (MF2)}} \\
\midrule
& \textbf{Species} & \textbf{Class} & \textbf{Prop. (\%)} &
  \textbf{Species} & \textbf{Class} & \textbf{Prop. (\%)} \\
\midrule
\multirow{5}{*}{\rotatebox{90}{\parbox{72pt}{\centering\textbf{Exploratory}}}}
& HMN       & Iso-alkane        & 26.2 & NBCH      & Cyclo-alkane      & 79.2 \\
& NC11H24   & N-alkane          & 12.3 & NC10H22   & N-alkane          & 3.6 \\
& NC15H32   & N-alkane          & 15.7 & NC11H24   & N-alkane          & 2.5 \\
& XC12H26   & Iso-alkane        & 21.1 & CPN       & Aldehyde/Ketone   & 10.7 \\
& CPN       & Aldehyde/Ketone   & 24.7 & C10H22-2 & Iso-alkane        & 4.0 \\
\midrule
\multirow{5}{*}{\rotatebox{90}{\parbox{72pt}{\centering\textbf{Practical}}}}
& O-XYL    & Aromatic          & 14.8 & NBCH      & Cyclo-alkane      & 77.6 \\
& NC13H28   & N-alkane          & 26.6 & NC12H26   & N-alkane          & 5.7 \\
& HMN       & Iso-alkane        & 26.3 & DECALIN   & Cyclo-alkane      & 9.0 \\
& NBCH      & Cyclo-alkane      & 11.5 & A2CH3     & Aromatic          & 2.0 \\
& C10H22-2 & Iso-alkane        & 20.8 & C11H24-2 & Iso-alkane        & 5.7 \\
\bottomrule
\end{tabular}
\end{adjustbox}
\vspace{-0.4cm}
\end{table}

\begin{figure}[t]
    \centering
    \includegraphics[width=1\linewidth]{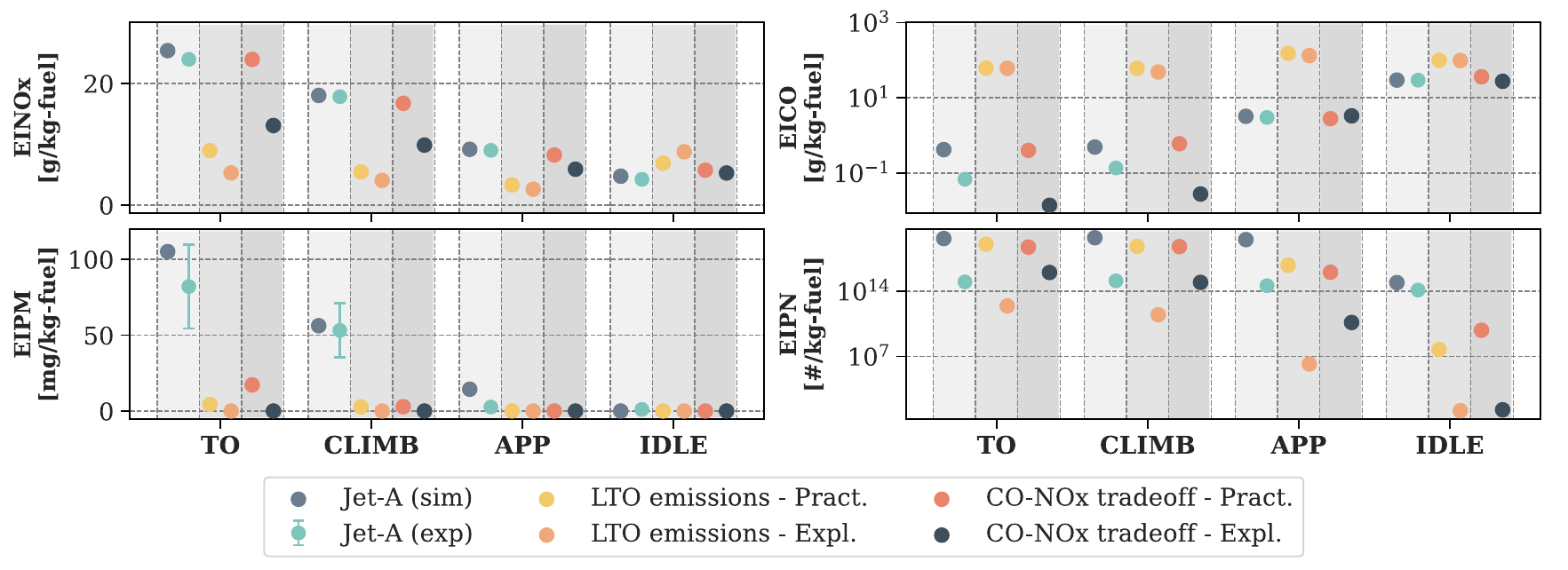}
    \scriptsize
    \caption{Results for the optimal fuel candidates identified by the \textit{Fuel Optimizer} for MF1 (CO--NO\textsubscript{x} trade-off) and MF2 (LTO regulated emissions), evaluated under both exploratory and practical scenarios. The selected blends were simulated using the 1D reactor network under LTO conditions and compared against both experimental and simulated Jet A data, which serves as the baseline reference for performance improvement.}
    \label{Reactor_network_optimum}
    \vspace{-0.5cm}
\end{figure}

\vspace{-0.4cm}
The optimum blends obtained for both MF1 and MF2 and following the exploratory and the practical approaches (which composition is described in table \ref{tab:fuel_blends}, see Supplementary Material for a glossary of species identifiers) have been simulated in the 1D reactor network model described in section \ref{subsec:hifi} at the four operating points of the LTO cycle. Fig.~\ref{Reactor_network_optimum} shows the reactor network results compared against both experimental measurements and simulations with Jet A. The results corroborate the findings of the \textit{Fuel Optimizer}. Braking the CO-NO$_x$ tradeoff (MF1) is challenging, with the optimum fuel blends for MF1 showing merit function values of -1.185 and -0.357 for the exploratory and the practical approaches, respectively. The optimum fuel of the exploratory approach shows lower CO emissions than Jet A at all conditions, and lower NO$_x$ at all conditions except idle (due to the positive effect of cyclopentanone of flame speed and flame temperature). However, the optimum fuel of the practical approach marginally improves NO$_x$ while showing similar CO levels than Jet A, suggesting that conventional hydrocarbons cannot break the CO-NO$_x$ tradeoff as oxygenated fuels do. 

The optimum fuel blends for MF2 showed almost zero soot mass and lower PN than Jet A at all conditions, suggesting that synthetic fuels mainly composed by alkanes shift the particle size distribution towards smaller diameters. While the optimum MF2 blend from the exploratory approach led to significantly lower PN values than Jet A, the optimum MF2 blend from the practical approach showed less benefits due to the aromatic content needed to meet the seal swell requirement. Interestingly, both optimum MF2 fuel candidates showed much higher CO emissions than Jet A. This is because CO was not included in MF2 and, therefore, CO emission values were not taken into consideration by the genetic algorithm. This result highlights the importance of a careful definition of the merit function and constraints, since the GA will exploit the gaps in the merit function to take advantage and provide good results. For this example, an additional constraint limiting the maximum allowable CO emissions will limit the ability of the GA to exploit the existing gap in MF2. However this is out of the scope of this paper, which main focus is to describe the \textit{Fuel Optimizer} showing two examples of its capability rather than propose the \textit{perfect} merit function for SATF design (which would be very subjective and user-specific).  

%%%%%%%%%%%%%%%%%%%%%%%%%%%%%%%%%%%%%%%%%%%%%%%%%%%%%%%%%%%%
\vspace{-0.4cm}
\subsection{Implications for SATF Certification}
\label{subsec:implications}
\vspace{-0.4cm}
The ASTM D4054 qualification process represents a significant financial and temporal bottleneck, requiring years of effort, millions of dollars, and thousands of gallons of fuel to advance a novel pathway to approval~\cite{Rumizen2021}. The inverse design methodology presented here effectively introduces a "Tier 0" computational gate, evaluating thousands of potential blend compositions in seconds to eliminate candidates that would predictably fail D1655 limits or deliver negligible emissions reductions, ensuring only the most promising formulations proceed to costly physical testing.

The optimization results also offer direct chemical guidance for SATF producers. Paraffinic fuels from HEFA and Fischer-Tropsch pathways offer thermal stability and low soot tendency but fail to meet drop-in seal swelling requirements due to their lack of aromaticity~\cite{Corporan2011}. Conventional aromatics address this but drive nvPM formation via the HACA mechanism~\cite{Frenklach2002}. The present results indicate that cycloalkanes, particularly n-butyl cyclohexane, and oxygenated species such as ketones represent a more favorable alternative, achieving simultaneous reductions in NO\textsubscript{x}, particulate mass, and particulate number without compromising operability. Synthetic pathways capable of selectively producing these chemical classes should therefore be prioritized.

%%%%%%%%%%%%%%%%%%%%%%%%%%%%%%%%%%%%%%%%%%%%%%%%%%%%%%%%%%%%
\vspace{-0.4cm}
\subsection{Limitations}
\label{sec:limitations}
\vspace{-0.4cm}
The accuracy of the 1D reactor model used to simulate engine emissions depends on specific assumptions about the chemical mechanisms and boundary conditions imposed to mimic the operation of a CFD56-7B27 turbofan (discussed in detail in~\cite{lopez-pintor_fuel_2026}). However, the model has been extensively validated against experimental data (see supplementary material,~\cite{lopez-pintor_fuel_2026} and~\cite{mohammed_modeling_2026}). In addition, simulations were performed for cruise-like operation (30$\%$ thrust), and fuel effects on performance may change at other operating conditions. However, cruise represents the majority of engine duty time and fuel trends at cruise can generally be extrapolated to other conditions~\cite{2006-01-1987}.
Finally, properties of multi-component fuel candidates were estimated using blending rules that have an inherent uncertainty. However, fuel property uncertainties only affect the ability of fuel candidates to meet the ASTM D1655 constrains, but they will not directly impact the emission results and, therefore, the merit function results shown in this paper.  
Future work to address these potential limitations includes: (1) integrating more sophisticated, non-ideal property blending rules~\cite{larranaga_data-driven_2026}, (2) expanding the emissions database to different engine operating conditions, and (3) validating selected fuel candidates by comparing their performance against Jet A in laboratory-scale facilities.  

%%%%%%%%%%%%%%%%%%%%%%%%%%%%%%%%%%%%%%%%%%%%%%%%%%%%%%%%%%%%
\vspace{-0.4cm}
{\section{Conclusions}
\label{sec:conclusions}
\vspace{-0.4cm}
This paper presents the \textit{Fuel Optimizer}, a data-driven machine learning-based numerical framework for designing, screening and optimization of synthetic aviation turbine fuel (SATF) formulations. The model is built on a palette of chemical components combined to generate a database of SATF-like blends that meet the ASTM D1655 standards for density, viscosity, distillation, freeze point, flash point and neat heat of combustion. Fuel blends of the database have been simulated in a 1D reactor network model of the CFM56-7B27 engine at cruise-like conditions (30$\%$ thrust) to generate a dataset of pollutant emissions and performance results that has been used to train a machine-learning surrogate model of the engine. The machine learning model demonstrated to be able to predict the 1D model results within 5$\%$ in most cases but at a fraction of the computational cost, allowing the evaluation of thousands of  fuel candidates per minute. A genetic algorithm was coupled with the machine learning model to optimize the composition of fuel candidates as per a user-defined merit function and while imposing a series of user-defined constraints. For the current version of the model, merit functions may include engine emissions, fuel properties and/or chemical species or hydrocarbon families. However, other parameters such as fuel cost or engine performance can and will be added in future work. 

The ability of the \textit{Fuel Optimizer} to explore potential SATF candidates was explored for two different merit functions: MF1, which was designed to break the CO-NO$_x$ trade-off typical of jet engines, and MF2, which was designed to minimize LTO-regulated emissions. For each merit function, two different set of constraints were imposed: an \textit{exploratory} set of constraints that included ASTM D1655 limits for density, viscosity, distillation, freeze point, flash point and neat heat of combustion but allowed any hydrocarbon class, and a \textit{practical} set of constraints that limited the amount of olefins to 5$\%$, did not allow oxygenated species, and required a minimum seal swelling rate equal to that of conventional kerosene. 

The \textit{Fuel Optimizer} was able to improve the merit function results of the database used for training and highly improve the merit function results of Jet A. A detailed analysis of the fuel candidates proposed by the genetic algorithm revealed that the CO-NO$_x$ trade-off can be broken by combining oxygenated species or species with a higher air-fuel stoichiometric ratio (lower flame temperature and lower NO$_x$) and species with high flame speed (less residence time required to complete combustion), with cyclopentanone being a key species for simultaneous reduction of CO and NO$_x$. Analyses also revealed that n-butylcyclohexane may act as a super-component to reduce LTO-regulated emissions (NO$_x$, PM and PN) when included in large amounts (> 70$\%$). Finally, results of the genetic algorithm indicate that small amounts of mono-aromatics (< $5\%$) combined with cyclo-alkanes is the preferred approach to meet seal swell constraints while minimizing soot emissions. Simulations of the optimum blends at LTO conditions in the 1D reactor network model supported the results from the \textit{Fuel Optimizer}, giving confidence in the ability of the \textit{Fuel Optimizer} to be used as a SATF formulation, screening and optimization tool.      

%%%%%%%%%%%%%%%%%%%%%%%%%%%%%%%%%%%%%%%%%%%%%%%%%%%%%%%%%%%%
\vspace{-0.4cm}
\section*{Data availability}
\vspace{-0.4cm}
The database used in this work is included in the paper supplementary material. The reactor network model will be available upon request. The machine learning surrogate and fuel optimizer are available at \url{https://github.com/alarran/fuel-optimizer}.

%%%%%%%%%%%%%%%%%%%%%%%%%%%%%%%%%%%%%%%%%%%%%%%%%%%%%%%%%%%%
\vspace{-0.4cm}
\section*{Acknowledgements}
\vspace{-0.4cm}
Sandia National Laboratories is a multi-mission laboratory managed and operated by National Technology $\&$ Engineering Solutions of Sandia, LLC (NTESS), a wholly owned subsidiary of Honeywell International Inc., for the U.S. Department of Energy’s National Nuclear Security Administration (DOE/NNSA) under contract DE-NA0003525. This written work is authored by an employee of NTESS. The employee, not NTESS, owns the right, title and interest in and to the written work and is responsible for its contents. Any subjective views or opinions that might be expressed in the written work do not necessarily represent the views of the U.S. Government. The publisher acknowledges that the U.S. Government retains a non-exclusive, paid-up, irrevocable, world-wide license to publish or reproduce the published form of this written work or allow others to do so, for U.S. Government purposes. The DOE will provide public access to results of federally sponsored research in accordance with the DOE Public Access Plan. ALJ and SLB acknowledge support from the National Science Foundation AI Institute in Dynamic Systems (grant number 2112085)
\vspace{-0.4cm}
%%%%%%%%%%%%%%%%%%%%%%%%%%%%%%%%%%%%%%%%%%%
%%% Bibliography
%%%%%%%%%%%%%%%%%%%%%%%%%%%%%%%%%%%%%%%%%%%

\bibliographystyle{elsarticle-num}
\end{document}

%%%%%%%%%%%%%%%%%%%%%%%%%%%%%%%%%%%%%%%%%%%%%%%%%%%%%%%%%%%%%
\appendix
\label{appendix}

\begin{center}
    \section*{\bfseries\fontsize{18}{20}\selectfont Supplementary Information}
    \vspace{-5mm}
    {\fontsize{15}{18}\selectfont The \textit{Fuel Optimizer}: A Data-Driven Numerical Framework for Formulation of  Aviation Turbine Fuel}

    Ana Larrañaga$^{1,2,*}$
    Steven L. Brunton $^{1,2}$
    Jacobo Porteiro$^3$
    Dario Lopez-Pintor$^4$
    
    \small 
    $^1$ Department of Mechanical Engineering, University of Washington, Seattle, WA 98195, United States \\
    $^2$ NSF AI Institute in Dynamic Systems, University of Washington, Seattle, WA 98195, United States \\
    $^3$ CINTECX, Universidade de Vigo, Grupo de Tecnoloxía Enerxética (GTE), Vigo, 36310, Spain \\
    $^4$ Sandia National Laboratories, 7011 East Ave, Livermore, 94550, CA, United States \\

\end{center}
\tableofcontents
\renewcommand{\thesection}{\arabic{section}}
\renewcommand{\thesubsection}{\arabic{section}\Alph{subsection}}
\vspace{-0.4cm}

\footnotetext[1]{Corresponding author: alarra@uw.edu}
%%%%%%%%%%%%%%%%%%%%%%%%%%%%%%%%%%%%%%%%%%%%%%%%%%%%%%%%%%%%%%%%%%%%%
\newpage
\section{1-D Reaction Network Validation} 
\label{SI:reaction-network-validation}
\begin{figure}[b]
    \centering
    \vspace{-0.5cm}
    \includegraphics[width=0.43\linewidth]{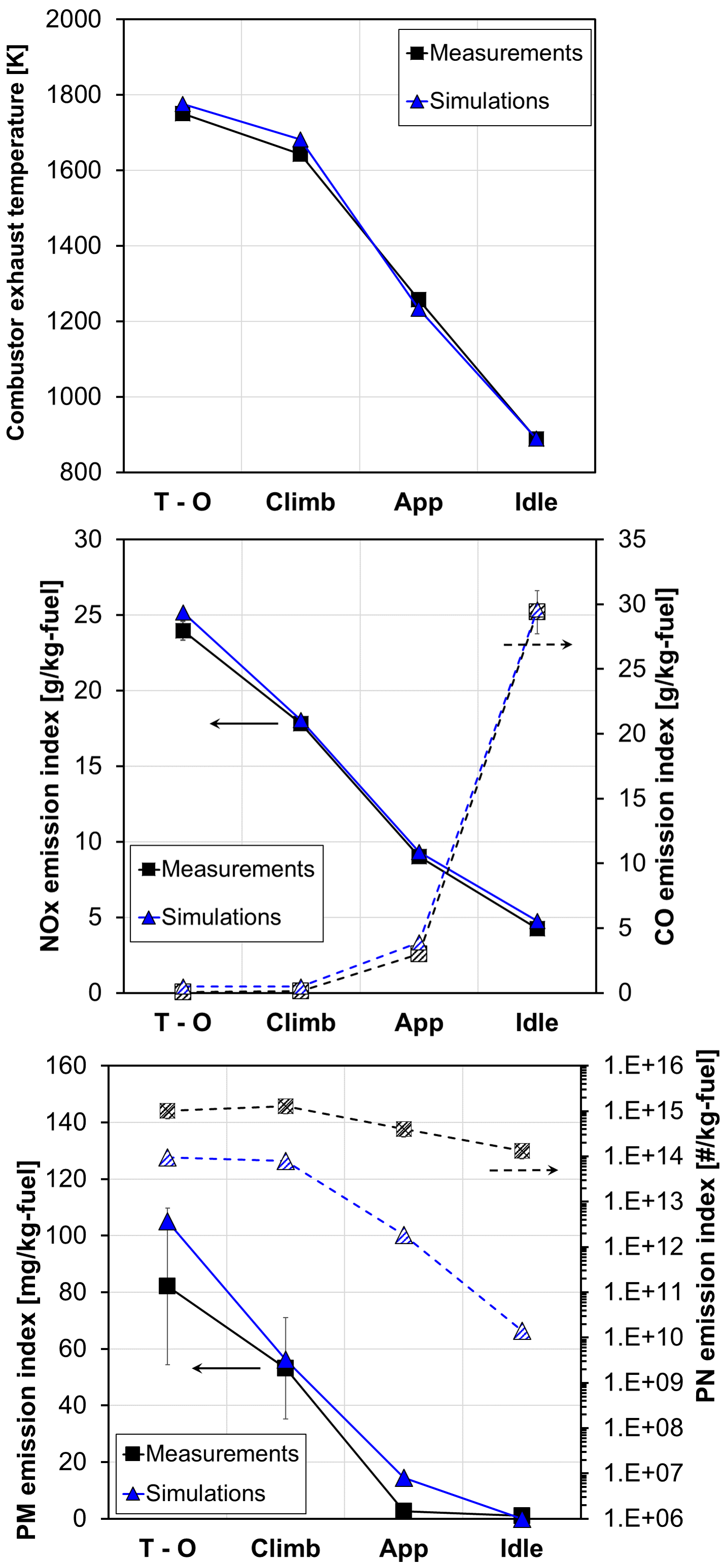}
    \scriptsize
    \vspace{-0.1cm}
    \caption{Comparison between reactor engine model simulations and experimental measurements of a CFM56-7B27 engine at LTO conditions with Jet A. Figure includes combustor outlet temperature (top), NO$_x$ and CO emissions (middle), and particulate mass and particle number emissions (bottom)}
    \label{figAppB}
\end{figure}

\vspace{-0.4cm}
The ability of the reactor network model to reproduce the emissions of the CFM56-7B27 engine at landing and take-off (LTO) conditions was validated by comparison against measurements from the ICAO Emissions databank~\cite{noauthor_icao_nodate}. The model was solved for a Jet A surrogate fuel~\cite{lopez-pintor_fuel_2026} at the four operating modes of the LTO cycle: take-off (100$\%$ thrust), climb (85$\%$ thrust), approach (30$\%$ thurst) and idle (7$\%$ thrust), and the results are shown together with those of the corresponding engine tests in fig. \ref{figAppB}. The reactor network was able to simulate combustor outlet temperature and gas-phase emissions (NO$_x$ and CO) with high accuracy. Particulate mass was also predicted within the uncertainty of the measurements. However, larger discrepancies were observed in particle number, with the overall trend well captured but the absolute values underpredicted by the model. Using jet A as a baseline fuel, the model is deemed to be sufficiently accurate to study fuel effects on emissions and perform fuel-to-fuel comparisons. 

%%%%%%%%%%%%%%%%%%%%%%%%%%%%%%%%%%%%%%%%%%%%%%%%%%%%%%%%%%%%%%%%%%%%%
\section{Surrogate Accuracy for Exploratory and Practical Scenarios}
\label{SI:surrogate-accuracy}
\begin{figure}[b]
     \centering
     \vspace{-0.5cm}
     \includegraphics[width=0.70\linewidth]{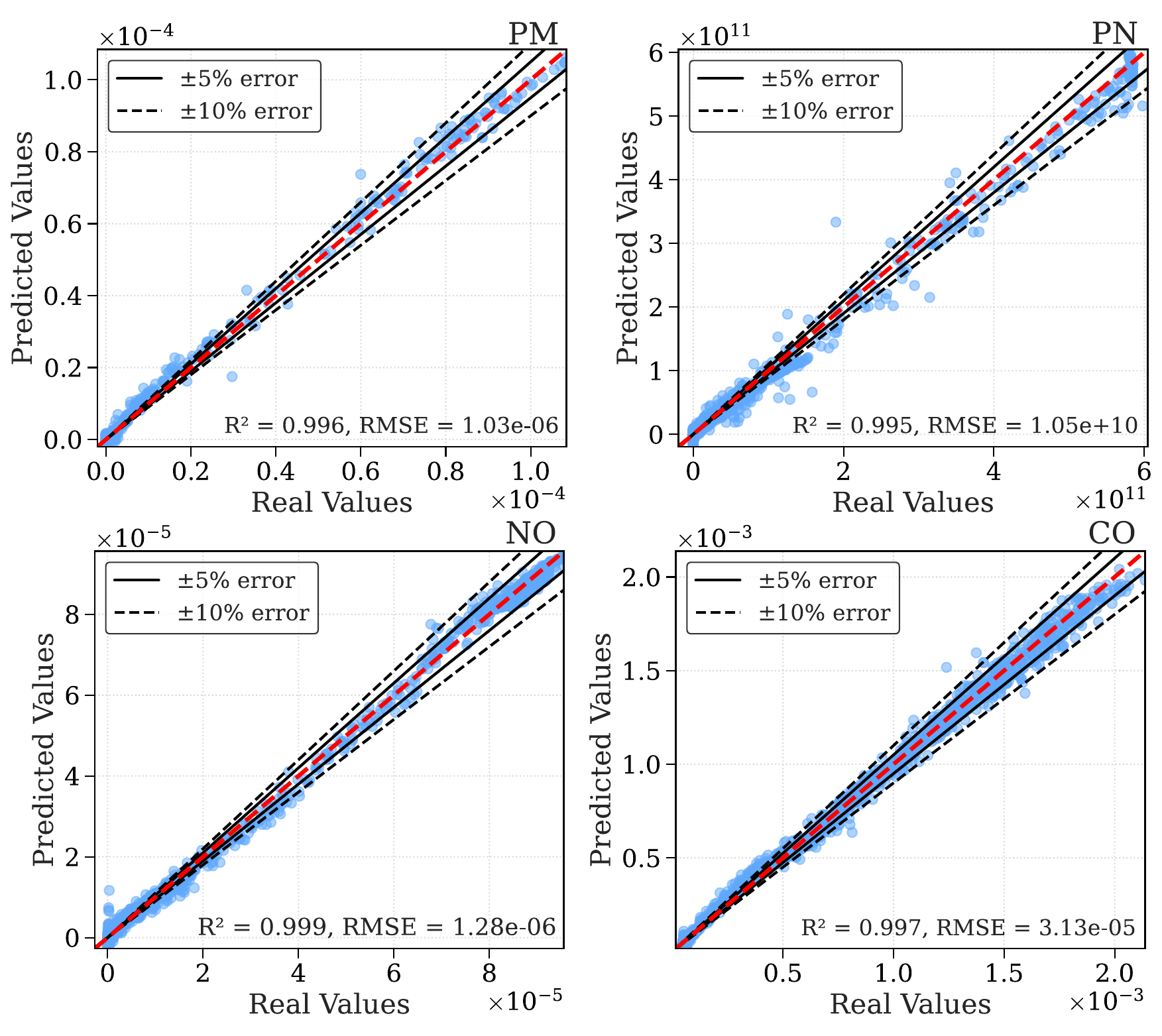}
     \scriptsize
     \vspace{-0.3cm}
     \caption{Surrogate model accuracy for the exploratory approach, where the database was not limited. The results are focused on the target variables relevant to the merit functions that are tested in this work.}
     \label{surrogate_performance}
 \end{figure}

\begin{figure}[t]
     \centering
     \includegraphics[width=0.70\linewidth]{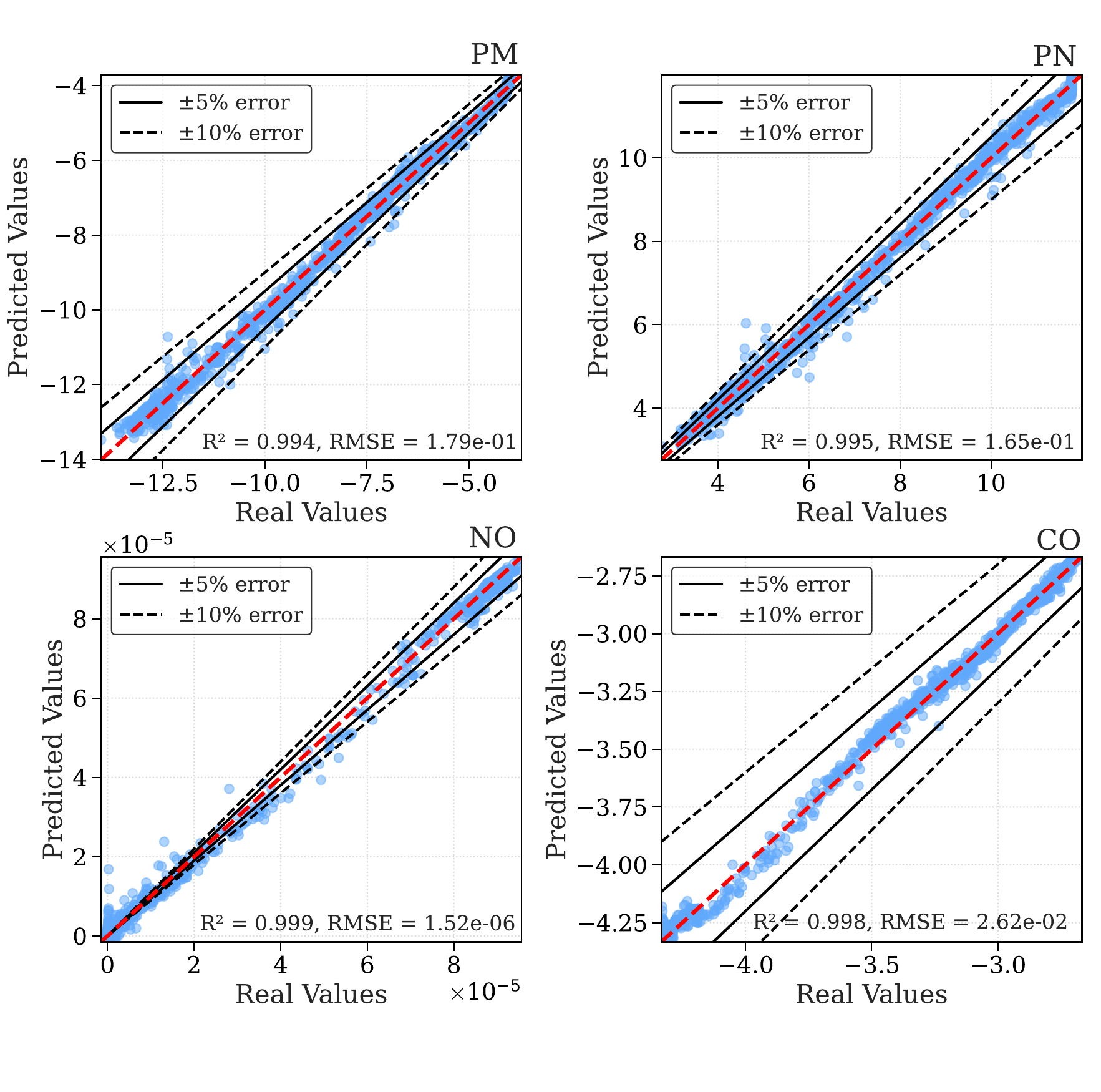}
     \scriptsize
     \vspace{-0.7cm}
     \caption{Surrogate model accuracy for the practical approach, where oxygenates were eliminated from the database. The results are focused on the target variables relevant to the merit functions that are tested in this work.}
     \label{surrogate_performance2}
     \vspace{-0.5cm}
 \end{figure}

\vspace{-0.4cm}
%%%%%%%%%%%%%%%%%%%%%%%%%%%%%%%%%%%%%%%%%%%%%%%%%%%%%%%%%%%%%%%%%%%%%
Two surrogate models were trained independently for the two optimization approaches. In both cases, the input normalization procedure and model architecture were optimized to minimize the prediction error. Separate models were required because the fuel databases differ between the two approaches: the exploratory approach includes oxygenated species, whereas the practical approach excludes them and imposes additional constraints on the allowable olefin content and on an additional fuel property. Figs.~\ref{surrogate_performance} and \ref{surrogate_performance2} show the performance of the jet engine surrogate model for the exploratory and the practical approaches, respectively. The models achieve prediction errors below 10\% for the targets relevant to MF1 and MF2, which include PM, PN, NO and CO. 

\textbf{Exploratory Fuel Optimizer}. The surrogate model is a fully connected feed-forward neural network mapping fuel blend composition to combustion emissions. The 51 inputs are species mole fractions (species absent from the database were excluded), min--max scaled. The architecture was selected via Bayesian hyperparameter optimization (Optuna, TPE sampler, 100 trials, 5-fold cross-validation) over number of layers (2--5), neurons per layer (16--128), activation function (\textsc{relu}, Leaky \textsc{relu}, \textsc{elu}), dropout rate (0.10--0.30), L1L2 coefficients ($10^{-6}$--$10^{-3}$), optimizer, and learning rate. The cross-validated optimum has three hidden layers of 56, 64, and 104 \textsc{relu} neurons, each followed by batch normalization (momentum $= 0.99$) and dropout (0.100, 0.109, 0.183), with L1L2 regularization ($l_1 = 1.15 \times 10^{-5}$, $l_2 = 6.0 \times 10^{-4}$) on all dense layers. The output layer has nine linear neurons predicting mole fractions of \ce{H2O}, \ce{NO}, \ce{NO2}, \ce{CO}, and \ce{CO2}; particle mass flow rate; average particle diameter; combustion temperature; and total particle number density (all min--max scaled). Training used RMSprop with MSE loss, batch size 24, up to 10{,}000 epochs, early stopping (patience $= 140$), and learning rate reduction on plateau (factor $= 0.1$, patience $= 15$, minimum LR $= 10^{-7}$). This configuration achieved the lowest mean cross-validation loss among the 100 trials and was retained as the final surrogate.

\textbf{Practical Fuel Optimizer}. This surrogate follows the same training framework as the exploratory model (Optuna TPE, 100 trials, 5-fold cross-validation, identical callbacks), differing in three aspects. First, the input layer receives 40 features, reflecting the reduced species space of the no-oxygenates blend database. Second, inputs are standardized (zero mean, unit variance) rather than min--max scaled. Third, three heavily skewed outputs --- \ce{CO} mole fraction, particle mass flow rate, and total particle number density --- are log$_{10}$-transformed before standardization; the remaining six targets are left on the original scale. The cross-validated optimum has two hidden layers of 120 and 128 \textsc{relu} neurons, each followed by batch normalization (momentum $= 0.99$) and dropout (0.289 and 0.136), with L1L2 regularization ($l_1 = 1.07 \times 10^{-4}$, $l_2 = 1.92 \times 10^{-5}$) on the hidden layers; the linear output layer uses fixed regularization ($l_1 = l_2 = 10^{-5}$). Training used RMSprop with MSE loss, batch size 24, early stopping (patience $= 140$), and learning rate reduction on plateau (factor $= 0.1$, patience $= 15$, minimum LR $= 10^{-7}$); the initial learning rate was tuned by Optuna over $10^{-5}$--$10^{-2}$. This configuration achieved the lowest mean cross-validation loss among the 100 trials and was retained as the final surrogate.

Details of the implementation can be found in \url{https://github.com/alarran/fuel-optimizer/tree/main/training-models} and the analysis of the results to generate the complete figure for all the outputs can be found in \url{https://github.com/alarran/fuel-optimizer/blob/main/AnalysisModel.ipynb}.

 %%%%%%%%%%%%%%%%%%%%%%%%%%%%%%%%%%%%%%%%%%%%%%%%%%%%%%%%%%%%%%%%%%%%%
 \vspace{-0.5cm}
\section{Optimization Algorithm} 
\label{SI:fuel-database}
\vspace{-0.4cm}
The fuel blend optimizer is a single-objective Genetic Algorithm (GA) implemented in \texttt{pymoo}. Each candidate is a nine-dimensional continuous vector: four variables in [0, 1] decoded into integer species indices over the full component library, and five variables in [0.02, 1.0] giving the corresponding volumetric blend fractions. The first blend component is fixed for an entire run, locking the search to a specific species from the fuel database; one GA run is executed per candidate species, so the procedure maps high-quality blends across all components rather than converging to a single global optimum. After each decode step, proportions are normalized to sum to unity, with a minimum fraction of 0.02 enforced per component by iterative redistribution. The population is 100 individuals; crossover and per-gene mutation probabilities are both 0.9, with mutations drawn uniformly within each variable's bounds by a custom \texttt{EnhancedMutation} operator. Duplicate candidates are eliminated at each generation. Feasibility is enforced through seven fuel-specification constraints (distillation temperatures, flash point, density, freezing point, viscosity, and net heat of combustion), extended to nine in the practical scenario by adding seal swell index and olefin content bounds, plus a sum-to-unity constraint; infeasible candidates receive an objective penalty of $10^6$. Termination occurs after 200 generations (exploratory) or 300 generations (practical), with a fallback loop that re-runs 50-generation sub-optimizations until at least one feasible candidate is recorded per species. The objective function analyzed in the present work are presented in the Methods section of the paper.

Details of the implementation can be found in \url{https://github.com/alarran/fuel-optimizer/tree/main/fuel_optimizer/optimizer}.

%%%%%%%%%%%%%%%%%%%%%%%%%%%%%%%%%%%%%%%%%%%%%%%%%%%%%%%%%%%%%%%%%%%%%
\vspace{-0.5cm}
\section{Full Chemical Family Analysis for Exploratory Approach} 
\label{SI:fuel-database}
\vspace{-0.4cm}

\begin{figure}[b]
    \centering
    \includegraphics[width=1\linewidth]{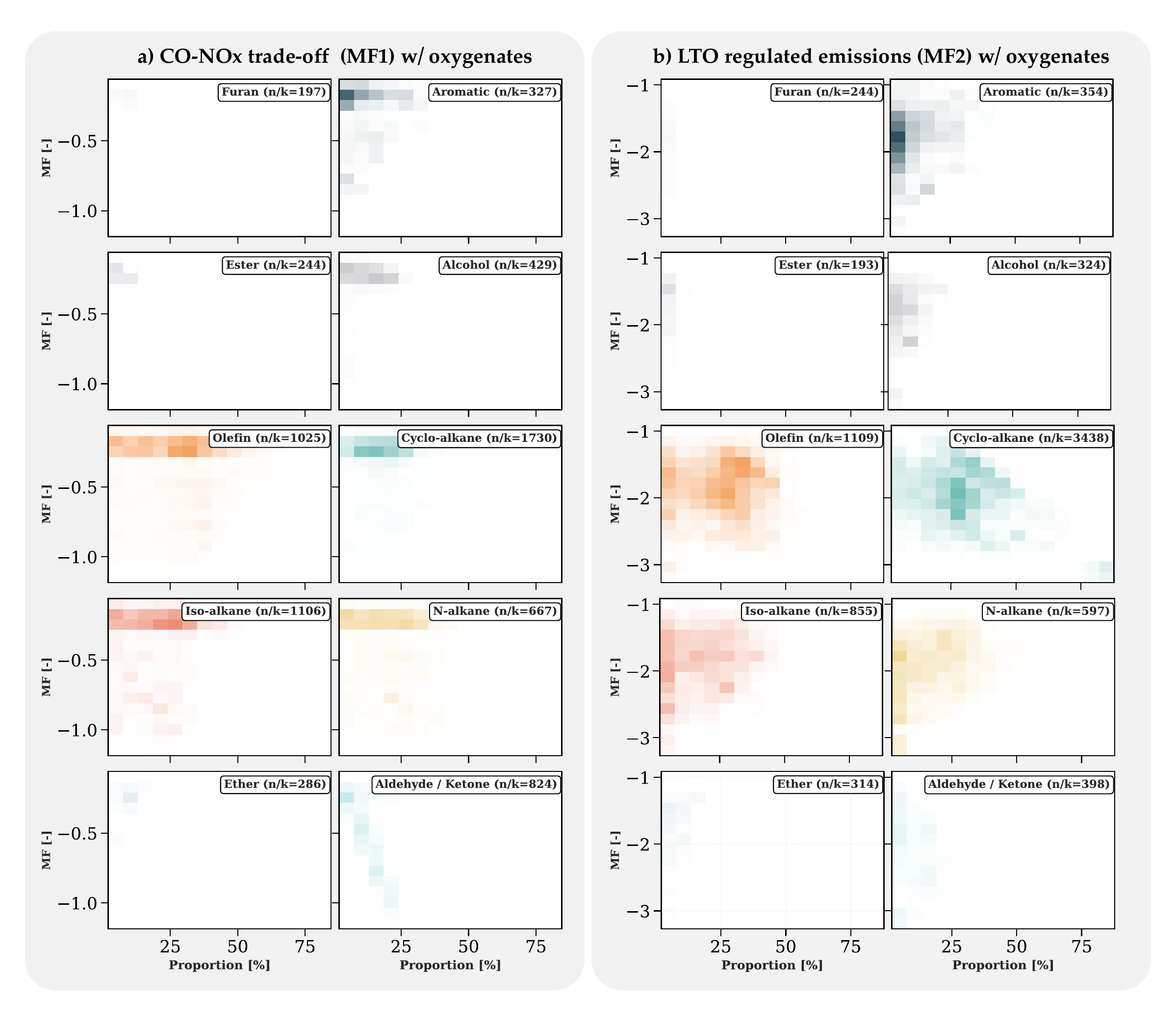}
    \scriptsize
    \vspace{-0.5cm}
    \caption{Full results for the \textit{exploratory approach} analyzed by chemical family.}
    \vspace{-0.5cm}
    \label{results-3-sup}
\end{figure}

Results for the \textit{Exploratory Fuel Optimizer} analyzed by chemical family. Fig.~\ref{results-3-sup} (a) on the left corresponds to the CO--NO\textsubscript{x} trade-off optimization, whereas figure (b) on the right corresponds to the LTO regulated-emissions optimization. The merit-function value is plotted on the y-axis, while the proportion of each chemical family in the blend is plotted on the x-axis. The color intensity represents the frequency of occurrence, with darker regions indicating a higher density of blends within a given composition range. The annotations in each subplot indicate the ratio between the number of occurrences of a given family among all candidate blends ($n$) and the number of species of that family available in the database ($k$), which serves as an indicator of the preference of the GA for each family.

%%%%%%%%%%%%%%%%%%%%%%%%%%%%%%%%%%%%%%%%%%%%%%%%%%%%%%%%%%%%%%%%%%%%%
\vspace{-0.5cm}
\section{Complete Single-Component Fuel Database} 
\label{SI:fuel-database}
\vspace{-0.4cm}
A palette of 51 single-component fuels, shown in table \ref{palette}, was used to create a database composed of 2508 multi-component blends that met the density, viscosity, distillation, freeze point, flash point and neat heat of combustion requirements for aviation fuel as defined in the ASTM 1655 standard~\cite{noauthor_standard_nodate}. Multi-component blend properties were estimated from the single-component fuel properties using blending rules, with the single-component fuel properties obtained from the National Institute of Standards and Technology (NIST) Standard Reference Database~\cite{noauthor_nisttrc_nodate}. More specifically, ideal liquid, Kendall-Monroe equation~\cite{kendall_viscosity_1917}, NIST REFPROF~\cite{lemmon_nist_2018}, freeze point blending index~\cite{coburn_determination_2022}, and Wickey-Chettenden equation~\cite{o_flash_1963} were used to calculate density and neat heat of combustion, viscosity, distillation, freeze point and flash point, respectively. 

The complete database used in this work is publicly accessible in: \url{https://github.com/alarran/fuel-optimizer}.  
  
\begin{table}[ht]
\centering
\scriptsize
\renewcommand{\arraystretch}{1.3}
\caption{Palette of single-component fuels used to generate the fuel database}
\label{palette}
\begin{tabularx}{\textwidth}{X|X|X|X}
\textbf{Hydrocarbon class} &
\textbf{Species name} &
\textbf{Identifier} &
\textbf{SMILES notation} \\
\hline
\hline
N-alkane                   & N-nonane              & NC9H20                         & CCCCCCCCC                          \\ \hline
N-alkane                   & N-decane              & NC10H22                         & CCCCCCCCCC                         \\ \hline
N-alkane                   & N-undecane             & NC11H24                         & CCCCCCCCCCC                        \\ \hline
N-alkane                   & N-dodecane              & NC12H26                        & CCCCCCCCCCCC                       \\ \hline
N-alkane                   & N-tridecane              & NC13H28                      & CCCCCCCCCCCCC                      \\ \hline
N-alkane                   & N-tetradecane            & NC14H30                      & CCCCCCCCCCCCCC                     \\ \hline
N-alkane                   & N-pentadecane          & NC15H32                        & CCCCCCCCCCCCCCC                    \\ \hline
N-alkane                   & N-hexadecane                & NC16H34                   & CCCCCCCCCCCCCCCC                   \\ \hline
N-alkane                   & N-octadecane                  & NC18H38                 & CCCCCCCCCCCCCCCCCC                 \\ \hline
N-alkane                   & N-eicosane                  & NC20H42                   & CCCCCCCCCCCCCCCCCCCC               \\ \hline
Iso-alkane                 & 2-methylheptane          & C8H18-2                     & CCCCCC(C)C                         \\ \hline
Iso-alkane                 & 2,2,5-trimethylhexane        & IC8                  & CC(C)CCC(C)(C)C                    \\ \hline
Iso-alkane                 & 2,2,4,6,6-pentamethylheptane  & XC12H26                  & CC(CC(C)(C)C)CC(C)(C)C             \\ \hline
Iso-alkane                 & 2,2,4,4,6,8,8-heptamethylnonane    & HMN            & CC(CC(C)(C)C)CC(C)(C)CC(C)(C)C     \\ \hline
Iso-alkane                 & 2-methyloctane          & C9H20-2                        & CCCCCCC(C)C                        \\ \hline
Iso-alkane                 & 2-methylnonane         & C10H22-2                         & CCCCCCCC(C)C                       \\ \hline
Iso-alkane                 & 2-methyldecane             & C11H24-2                     & CCCCCCCCC(C)C                      \\ \hline
Cyclo-alkane               & Butylcyclohexane          & NBCH                     & CCCCC1CCCCC1                       \\ \hline
Cyclo-alkane               & Decalin                 & DECALIN                       & C1CCC2CCCCC2C1                     \\ \hline
Aromatic                   & Toluene             & C6H5CH3                           & CC1=CC=CC=C1                       \\ \hline
Aromatic                   & Ethylbenzene            & C6H5C2H5                      & CCC1=CC=CC=C1                      \\ \hline
Aromatic                   & P-xylene            & P-XYL                           & CC1=CC=C(C=C1)C                    \\ \hline
Aromatic                   & O-xylene               & O-XYL                        & CC1=CC=CC=C1C                      \\ \hline
Aromatic                   & 1-ethyl-3-methylbenzene      & CH3C6H4C2H5                 & CCC1=CC=CC(=C1)C                   \\ \hline
Aromatic                   & 1,2,4 trimethylbenzene       & T124MBZ                  & CC1=CC(=C(C=C1)C)C                 \\ \hline
Aromatic                   & Tetralin                & TETRA                       & C1CCC2=CC=CC=C2C1                  \\ \hline
Aromatic                   & Naphthalene             & NAPH                       & C1=CC=C2C=CC=CC2=C1                \\ \hline
Aromatic                   & Alpha-methylnaphthalene       & A2CH3                 & CC1=CC=CC2=CC=CC=C12               \\ \hline
Aromatic                   & 1,2-diphenylethane          & C14H14                   & C1=CC=C(C=C1)CCC2=CC=CC=C2         \\ \hline
Olefin                     & 4-hexadecene             & C16H32-4                      & CCCCCCCCCCC/C=C\textbackslash{}CCC \\ \hline
Olefin                     & 4-octadecene        & C18H36-4                           & CCCCCCCCCCCCC/C=C/CCC              \\ \hline
Olefin                     & 4-eicosene         & C20H40-4                            & CCCCCCCCCCCCCCC/C=C/CCC            \\ \hline
Olefin                     & 2,4-dimethyl-1,3-pentadiene     & I24C7D13               & CC(=CC(=C)C)C                      \\ \hline
Olefin                     & 1-nonene            & C9H18-1                            & CCCCCCCC=C                         \\ \hline
Olefin                     & 1-butylcyclohexene       & NBCH-1ENE                      & CCCCC1=CCCCC1                      \\ \hline
Olefin                     & 1,2,3,4,4a,5,6,8a-octahydro-naphthalene     &  DCLD1  & C1CCC2=C(C1)CCCC2                  \\ \hline
Olefin                     & 2-methyl-1-nonene      & C10H20-1-2                        & CCCCCCCC(=C)C                      \\ \hline
Olefin                     & 2,4,4,6,6 pentamethyl-1-heptene   &  AC12H24            & CC(=C)CC(C)(C)CC(C)(C)C            \\ \hline
Olefin                     & 2,4,4,6,8,8 hexamethyl-1-nonene     &   AC15H30         & CCC(C)(C)C(C)(C)CC(=CC(C)C)C       \\ \hline
Olefin                     & 2,2,4,4,8,8 hexamethyl-6-methylenenonane  & BC16H32     & C=C(C)CC(C)CC(C)CC(C)CC(C)C        \\ \hline
Alcohol                    & Ethanol             & C2H5OH                           & CCO                                \\ \hline
Alcohol                    & Methanol                    & CH3OH                   & CO                                 \\ \hline
Alcohol                    & 1-pentanol              &    NC5H11OH                    & CCCCCO                             \\ \hline
Alcohol                    & 3-methyl-1-butanol    &    IC5H11OH                      & CC(C)CCO                           \\ \hline
Ether                      & Anisole                   &  C6H5OCH3                    & COC1=CC=CC=C1                      \\ \hline
Ester                      & Methylbutanoate        & MB                       & CCCC(=O)OC                         \\ \hline
Aldehyde / Ketone          & Cyclopentanone         & CPN                        & C1CCC(=O)C1                        \\ \hline
Aldehyde / Ketone          & Valeraldehyde         &     PC4H9CHO                     & CCCCC=O                            \\ \hline
Aldehyde / Ketone          & Hexanal            &   NC5H11CHO                          & CCCCCC=O                           \\ \hline
Aldehyde / Ketone          & 2,4-dimethyl-3-pentanone     &   NEC7Y4                & CC(C)C(=O)C(C)C                    \\ \hline
Furan                      & Furfural              &      F2CHO                    & C1=COC(=C1)C=O                     \\ \hline

\end{tabularx}
\end{table}

%%%%%%%%%%%%%%%%%%%%%%%%%%%%%%%%%%%%%%%%%%%
%%% Bibliography
%%%%%%%%%%%%%%%%%%%%%%%%%%%%%%%%%%%%%%%%%%%

\bibliographystyle{elsarticle-num}